\documentclass[12pt]{article}
\usepackage{epsfig}
\usepackage{color}
\usepackage{epstopdf}
\usepackage{graphics}
\usepackage{amsmath}
\usepackage{amssymb}
\usepackage{mathrsfs}

\newtheorem{theo}{Theorem}[section] 
\newtheorem{lemma}[theo]{Lemma} 
\newtheorem{mydef}[theo]{Definition} 
\newtheorem{cor}[theo]{Corollary} 
\newtheorem{prop}[theo]{Proposition} 
\newtheorem{rmk}[theo]{Remark}

\newtheorem{hypo}[theo]{Hypothesis}

\newcommand{\Sone}{{\bf S}^1}

\newcommand{\SEtwo}{{\bf SE}(2)}

\newcommand{\ep}{\varepsilon}

\newcommand{\proof}{\noindent {\bf Proof} \hspace{0.2in}} 
\newcommand{\qed}{\hfill\mbox{\raggedright\rule{.07in}{.1in}}
  \vspace{1ex}} 

\newcommand{\Section}[1]{\section{#1} \setcounter{equation}{0}}

\title{Dynamics of meandering spiral waves with weak lattice perturbations}
 
\author{Petko Kitanov\\Department of Mathematics and Statistics\\
University of Ottawa\\Ottawa, ON K1N 6N5\\CANADA  
\and Victor G. LeBlanc\\Department of
  Mathematics and Statistics\\University of Ottawa\\Ottawa, 
ON K1N 6N5\\CANADA}
\date{\today}

\begin{document}

\maketitle

\begin{abstract}
Re-entrant spiral waves are observed in many different situations in nature, perhaps most importantly in excitable electrophysiological tissue where they are 
believed to be responsible for pathological conditions such as cardiac arrhythmias, epileptic seizures and hallucinations. Mathematically, spiral waves occur as
solutions to systems of reaction-diffusion partial differential equations (RDPDEs) which are frequently used as models for electrophysiological phenomena. 
Because of the invariance of these RDPDEs with respect to the Euclidean group $\SEtwo$ of planar translations and rotations, 
much progress has been made in understanding the dynamics and bifurcations
of spiral waves using the theory of group-equivariant dynamical systems. In reality however, Euclidean symmetry is at best an approximation. Inhomogeneities and
anisotropy in the medium of propagation of the waves break the Euclidean symmetry, and can lead to such phenomena as anchoring and drifting.  In this paper, we study
the effects on quasi-periodic meandering spiral waves of a small perturbation which breaks the continuous $\SEtwo$ symmetry, but preserves the symmetry of a regular square
lattice. 

\end{abstract}

\pagebreak
\Section{Introduction}

Spiral waves have been studied now extensively for a few decades since the pioneering work of Winfree \cite{Winfree1,Winfree2}.  One of the most important reasons why they have been given 
so much attention is the mounting evidence that in electrophysiological tissue (such as the myocardium, the visual cortex or the brain), spiral waves are typically symptomatic of pathological
conditions such as fibrillation, hallucinations or epileptic seizures \cite{Davidenko,Huangetal,KeenerSneyd,Mesin,Pertsov,Viventi}. 

A large class of mathematical models in which spiral waves occur as solutions are reaction-diffusion partial differential equations (RDPDEs) on planar domains
\begin{equation}
\frac{\partial U}{\partial t}(x,y,t)=D\cdot \nabla^2 U(x,y,t)+f(U(x,y,t),\lambda),
\label{basicRDPDE}
\end{equation}
where $U:\mathbb{R}^2\times\mathbb{R}^{+}\longrightarrow \mathbb{R}^N$ may represent for example electric potentials of different ions through cellular membranes, or concentrations of chemicals. The $N\times N$ constant
matrix $D$ determines diffusion coefficients, and $\nabla^2=\frac{\partial^2}{\partial x^2}+\frac{\partial^2}{\partial y^2}$ is the Laplacian.  The smooth function $f:\mathbb{R}^N\times\mathbb{R}^m
\longrightarrow\mathbb{R}^N$ describes the local reaction kinetics of the model, and $\lambda\in\mathbb{R}^m$ are model parameters.

When (\ref{basicRDPDE}) is posed on the whole plane, this equation admits an important symmetry property\footnote{In fact, (\ref{basicRDPDE}) also admits reflectional symmetries, but these will not be relevant for the
purposes of this paper.}: whenever $U(x,y,t)$ is a solution of (\ref{basicRDPDE}), then so is
\[
S(x,y,t)=U(x\cos\,\phi-y\sin\,\phi + p_1,x\sin\,\phi+y\cos\,\phi+p_2,t)
\]
for any angle $\phi$ in the circle group $\Sone$, and any $(p_1,p_2)\in\mathbb{R}^2$.  In the jargon of dynamical systems theory, we say that the right-hand side of (\ref{basicRDPDE}) is {\em $\SEtwo$-equivariant}, where $\SEtwo$ designates the group of all planar
translations and rotations.

Barkley \cite{B92,B94,BK,BKT} was the first to recognize the importance of $\SEtwo$ symmetry in describing the dynamics and bifurcations of spiral waves which were observed in numerical simulations
of (\ref{basicRDPDE}) and in physical experiments \cite{Li}.  Essentially, Barkley argued that the continuous symmetries of $\SEtwo$ should lead to low-dimensional (finite) ordinary differential equations models
to describe the basic modes of propagation of spiral waves: uniform spatial rotation of the spiral wave around a fixed point in space, or two-frequency epicyclic meandering spiral waves, where the tip of the
spiral wave rotates around a point which itself precesses (see for example Figures \ref{grideffects} and \ref{fig2} below).  Later, Sandstede, Scheel and Wulff \cite{SSW1,SSW2} proved mathematically how Barkley's finite
dimensional models are, in fact, center manifold reductions of (\ref{basicRDPDE}) in the context of infinite-dimensional $\SEtwo$-equivariant dynamical systems.   See also \cite{GLM1,GLM2}.
In this context, rigidly rotating spiral waves are examples of {\em relative equilibria}, and two-frequency epicyclic meandering waves are examples of {\em relative periodic solutions} (see \cite{SSW1,SSW2} 
for the precise mathematical definitions of these concepts).

In any real physical situation, it is clear that $\SEtwo$ symmetry is, at most, an approximation.  Boundaries, non-planar geometry, inhomogeneities and anisotropy (which are all important in 
the electrophysiological contexts alluded to earlier) break translational and/or rotational symmetries.  Moreover, there are several experimental studies illustrating how broken
translational and/or rotational symmetries may lead to dynamic states for spiral waves which are inconsistent with Euclidean symmetry: anchoring of spiral waves around regions of
inhomogeneities \cite{Davidenko,MMV98}, drifting along boundaries \cite{YP86,ZM96}, and phase-locking and/or drifting of two-frequency meandering waves in the bidomain model of cardiac electrophysiology \cite{Roth1,Roth2}.  These experimental observations led to 
a program of research \cite{BLM,BEL,CL,LeBlanc,LeBlancRoth,LW} which sought to use the theoretical apparatus of the center-manifold theorems of \cite{SSW1,SSW2} to characterize generic effects of {\em forced symmetry-breaking} on 
relative equilibria and relative periodic solutions.  As a basic paradigm to illustrate forced symmetry-breaking from $\SEtwo$, consider the following perturbation of (\ref{basicRDPDE})
\begin{equation}
\frac{\partial U}{\partial t}(x,y,t)=D\cdot \nabla^2 U(x,y,t)+f(U(x,y,t),\lambda)+\ep g(U(x,y,t),x,y,\ep),
\label{basicRDPDEpert}
\end{equation}
where $\ep$ is a small parameter and $g$ is some bounded function which depends non-trivially and explicitly on the spatial coordinates $x$ and $y$.  Whereas (\ref{basicRDPDEpert}) admits
full $\SEtwo$ symmetry when $\ep=0$, it typically only admits a subgroup $\Sigma\subset\SEtwo$ when $\ep\neq 0$, however small $\ep$ may be.  Depending on the physical situation that
one wants to describe, the subgroup $\Sigma$ is prescribed in advance.  Using this approach, it was shown in \cite{BLM,LeBlanc,LW} that spiral anchoring, boundary drifting, and phase-locking and/or 
drifting of meandering
waves in anisotropic media are generic consequences of forced symmetry-breaking from $\SEtwo$ to $\Sigma$, for appropriate choices of $\Sigma$.  In some sense, we can view the 
dynamical system generated by (\ref{basicRDPDEpert}) with $\ep\neq 0$ as being
close to an $\SEtwo$-equivariant dynamical system, so we expect solutions of (\ref{basicRDPDEpert}) to retain some of the features of the $\SEtwo$ symmetry (at least on a transient level), as well as features of a $\Sigma$-equivariant dynamical system.

\subsection{Lattice symmetry-breaking}

\begin{figure}[htpb]
\begin{center}
\includegraphics[width=2.0in]{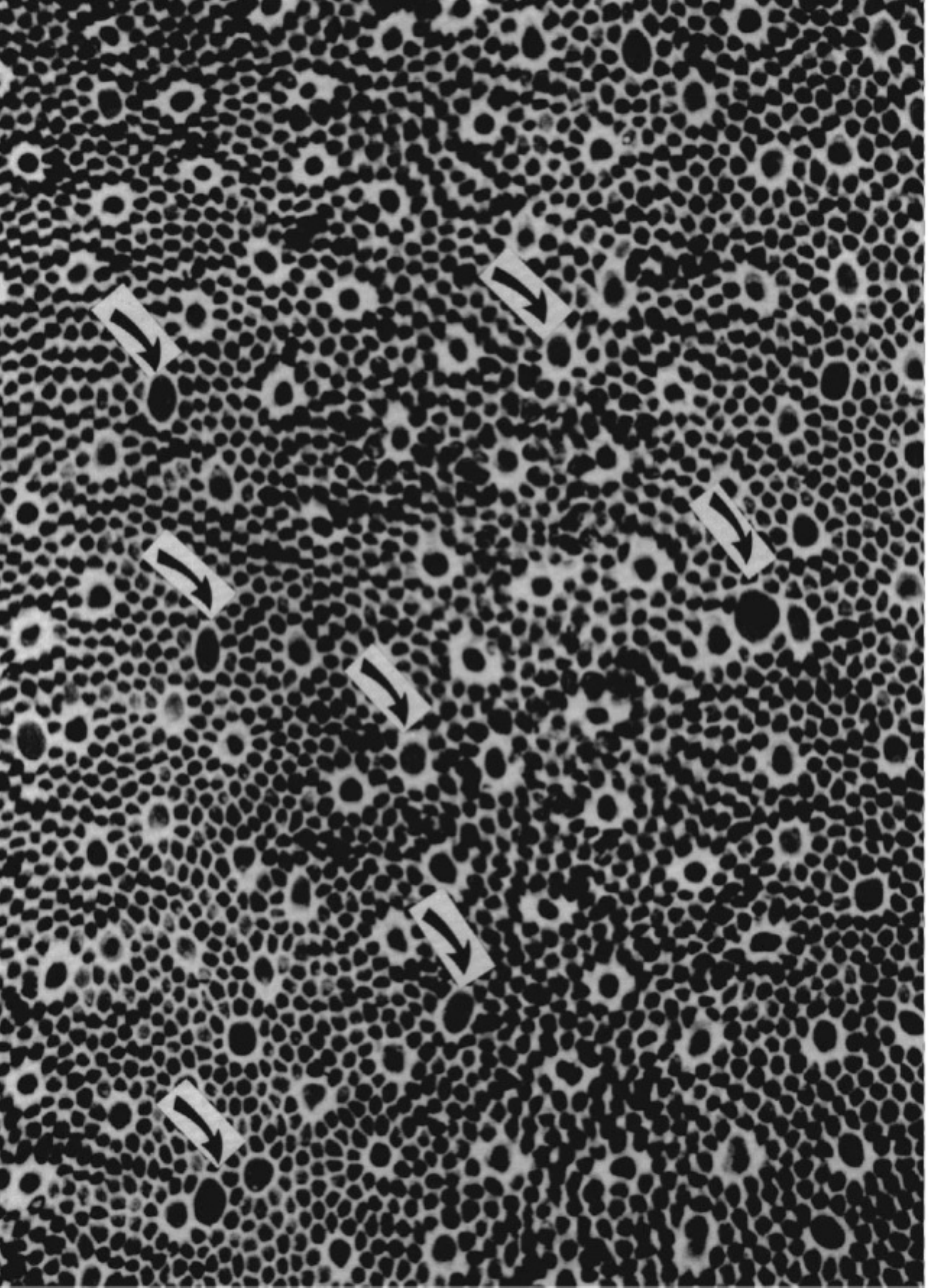}
\end{center}
\caption{Figure 2 of \cite{Ahnelt} (with permission) illustrating the spatial arrangement of cones in a cross section of the human retina.}
\label{retina}
\end{figure}

In a recent paper \cite{CL}, we studied the case where $\Sigma$ in (\ref{basicRDPDEpert}) is the group of rotational and translational symmetries of a regular square lattice.  Specifically, we
characterized the effects of this symmetry-breaking on relative equilibria (i.e. uniformly rigidly rotating spiral waves, or linearly translating waves with retracting tip).
The motivation was two-fold.  First, whenever an RDPDE model such as (\ref{basicRDPDE}) is used to describe electrophysiological waves, it is understood that this continuum model
is obtained after some averaging and homogenization over the cellular structure of the tissue.  However, if the size of the spiral wave (e.g. spiral step size) is comparable to the size of the cellular
structure, then inhomogeneities generated by gap junctions between cells may have an influence on the dynamics of the spiral wave.  
See for example Figure \ref{retina} (taken from \cite{Ahnelt}) which illustrates the spatial distribution of cones in a cross section of the human retina.
As a first approximation, one may assume that
the cells are arranged in a square grid, although one could argue that in certain regions of Figure \ref{retina}, perhaps a hexagonal grid is closer to reality.  In any event, our goal is to
characterize the generic first order effects of a lattice structure on spiral dynamics, so the specific geometry of the lattice (i.e. square vs hexagonal) is secondary in light of this. We choose a square
lattice to simplify the presentation, and for the second motivation described below.

A second motivation to study square lattice symmetry-breaking from $\SEtwo$ is to try to characterize effects of coarse spatial grids on the resolution of spiral wave dynamics in numerical simulations of (\ref{basicRDPDE}).  Although the RDPDE (\ref{basicRDPDE}) enjoys full $\SEtwo$ symmetry, numerical integration schemes which involve finite difference spatial discretization possess only
the symmetries of the lattice.  So it is not unreasonable to expect that numerical integrations of (\ref{basicRDPDE}) may reflect features which are 
characteristic of this grid symmetry, while still being ``close'' to being fully $\SEtwo$ symmetric \cite{BEL}.  This is especially true if the grid is coarse in relation
to the size of the spiral.  As an illustration of this point, in Figure \ref{grideffects} below, we show the path of the tip of a meandering spiral wave for four different numerical
integrations of the FitzHugh-Nagumo
system
\begin{equation}
\begin{array}{rcl}
{\displaystyle\frac{\partial u}{\partial t}}&=&{\displaystyle\nabla^2 u+\frac{1}{\tau}\left(u-\frac{1}{3}u^3-v\right)}\\[0.1in]
{\displaystyle\frac{\partial v}{\partial t}}&=&{\displaystyle \tau (u+\beta-\gamma v)}
\end{array}
\label{FNPDEup}
\end{equation}
on the domain $[-10\pi,10\pi]^2$ using a finite difference scheme (explicit in time) with varying grid sizes, and Neumann boundary conditions.  In particular, there is phase-locking of a 
four-petal closed meandering path for a coarse grid, and this is inconsistent with Euclidean symmetry \cite{BKT}.
\begin{figure}[htpb]
\begin{center}
\includegraphics[width=3.2in]{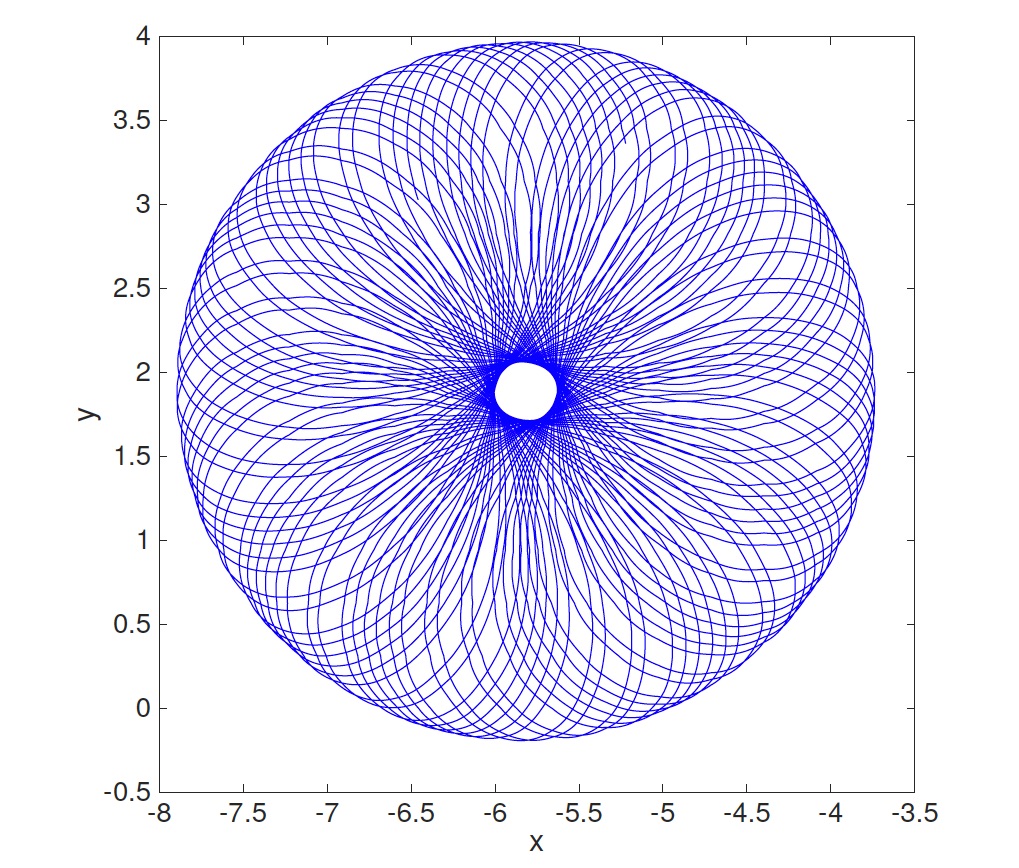}
\includegraphics[width=3.2in]{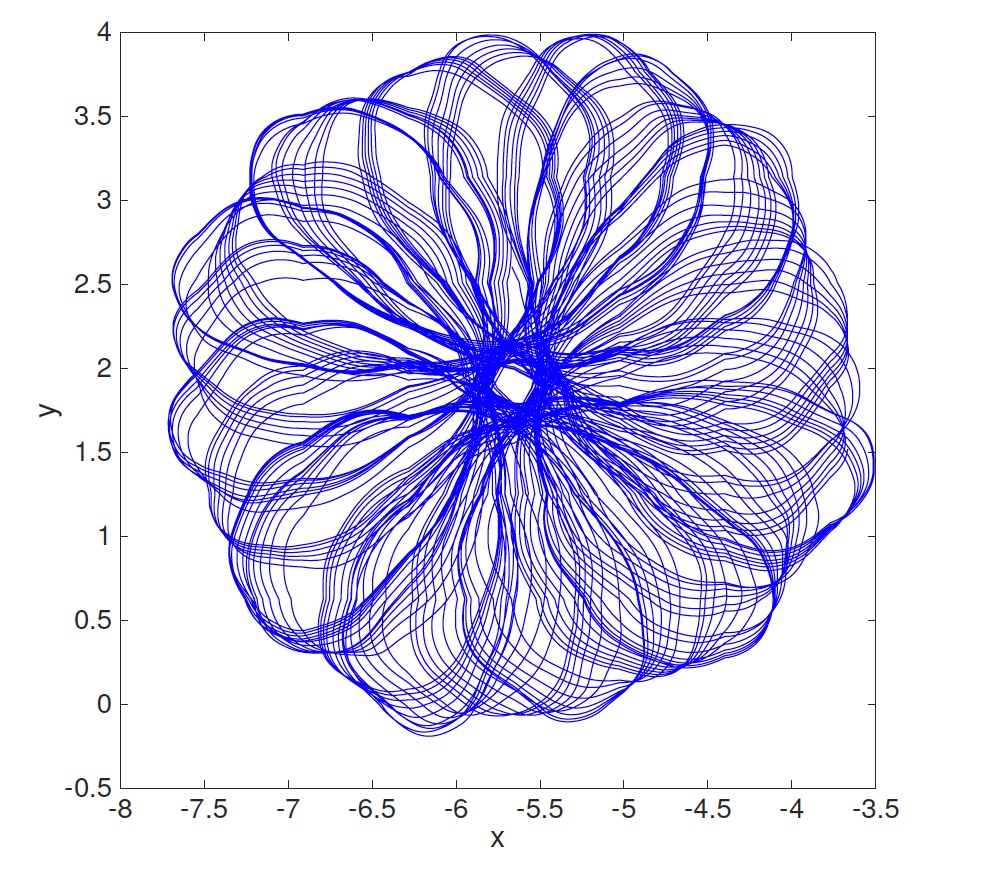}

\includegraphics[width=3.2in]{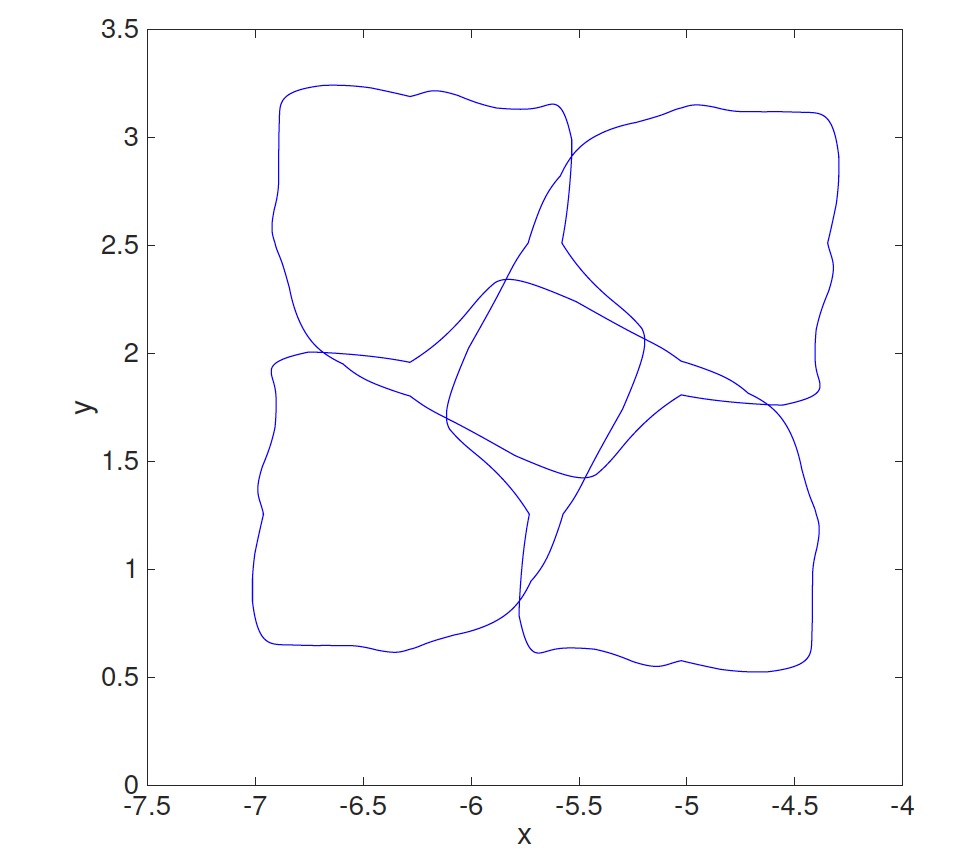}
\includegraphics[width=3.2in]{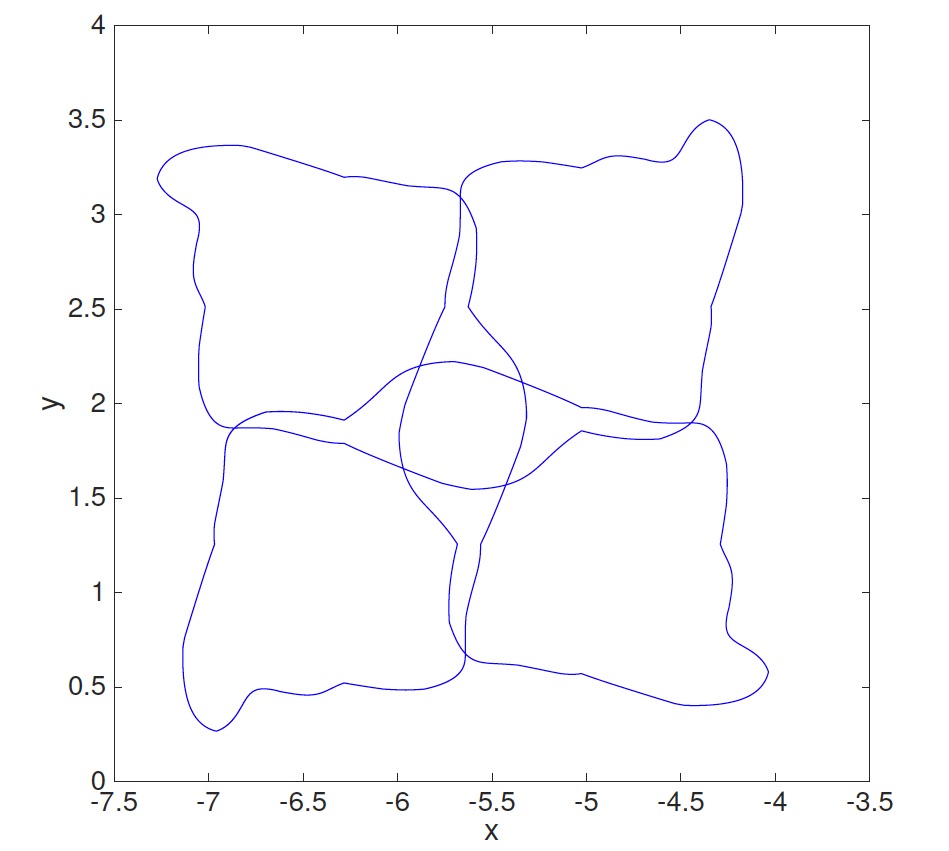}
\end{center}
\caption{Path of the spiral tip path for numerical integrations of (\ref{FNPDEup}) on the domain $[-10\pi,10\pi]^2$ using a finite difference scheme (explicit in time) with varying grid sizes:
$200\times 200$ (top left), $100\times 100$ (top right), $50\times 50$ (bottom two figures). The same initial condition is used in all cases, and transients have been removed.  The kinetic
parameters are $\beta=0.755$, $\gamma=0.5$ (for all four figures), and $\tau=0.2$ (both top figures and bottom left figure), or $\tau=0.25$ (bottom right). The bottom figures illustrate
a phase-locked four-petal closed meandering path with four-fold rotational symmetry, a result we will see in section 6 that is consistent with lattice symmetry-breaking.}
\label{grideffects}
\end{figure}
See also \cite{Adamatsky, Agladze, Sole, Xu} for further motivation as to the importance of understanding spiral wave dynamics in lattice structures.

As mentioned above, our emphasis in \cite{CL} was on characterizing the effects of lattice symmetry-breaking on relative equilibria. The goal of this present paper is to study the effects of this type of perturbation
on relative periodic solutions (otherwise known as modulated rotating waves, or modulated travelling waves \cite{B92,B94,BKT}).

\subsection{Numerical simulations}

At several places in this paper, we will present results of numerical simulations to illustrate some of the features which are predicted from our analysis.  
These simulations are numerical integrations of the perturbed FitzHugh-Nagumo system
\begin{equation}
\begin{array}{rcl}
{\displaystyle\frac{\partial u}{\partial t}}&=&{\displaystyle\nabla^2 u+\frac{1}{\tau}\left(u-\frac{1}{3}u^3-v\right)+\ep\,g_1(x,y)}\\[0.1in]
{\displaystyle\frac{\partial v}{\partial t}}&=&{\displaystyle \tau (u+\beta-\gamma v+\ep\,g_2(x,y))}
\end{array}
\label{FNPDE}
\end{equation}
where $\nabla^2=\frac{\partial^2}{\partial x^2}+\frac{\partial^2}{\partial y^2}$, $\tau$, $\beta$ and $\gamma$ are model parameters (which will be varied across the different simulations), the functions $g_1$ and $g_2$ are the lattice symmetry-breaking terms, and $\ep\geq 0$ is a small parameter.  When $\ep=0$, (\ref{FNPDE}) reduces to (\ref{FNPDEup}), and is symmetric under the planar group of rotations and translations, $\SEtwo$.  
We will consider functions $g_{1,2}$ of the form
\begin{equation}
\begin{array}{lll}
g_{i}=A_i+B_i (\cos(x/2)+\cos(y/2)) + C_i (\cos((3x-2y)/2)+\cos((2x+3y)/2)),\,\,\,i=1,2,
\end{array}
\label{gspecif}
\end{equation}
where the coefficients $A_i, B_i$ and $C_i$ will vary from simulation to simulation.  
Therefore, when $\ep\neq 0$, the $\SEtwo$ symmetry of (\ref{FNPDE}) is broken, but the rotational and translational symmetries of a square lattice are preserved.
Simulations with more harmonic components and/or with sine components were also performed, but did not lead to any dynamics significantly different from that obtained using 
(\ref{gspecif}).
\begin{mydef}
For the perturbation functions $g_i$ in (\ref{gspecif}), system (\ref{FNPDE}) is invariant under transformations which are combinations of
translations along the $x$ or the $y$ directions by integer multiples of $4\pi$, and rotations around the origin by $\pi/2$.  We will therefore refer to the points
\[
\{\,(4\pi n_1,4\pi n_2)\,\,|\,\,n_1,n_2\in\mathbb{Z}\,\}
\]
as {\em lattice points} in this paper.  System (\ref{FNPDE}) is also invariant under rotation by $\pi/2$ around the point $(2\pi,2\pi)$.  Hence, we will refer to the points
\[
\{\,(2\pi+4\pi n_1,2\pi+4\pi n_2)\,\,|\,\,n_1,n_2\in\mathbb{Z}\,\}
\]
as {\em dual lattice points} in this paper.
\end{mydef}
Throughout this paper, in the several figures where we present results of numerical simulations of (\ref{FNPDE}), we superimpose on the figure red squares and black crosses to guide the eye
as to the position of the lattice points and the dual lattice points respectively (see Figure \ref{fig4} for example).

For the numerical simulations which will be illustrated in this paper, the partial differential equation (\ref{FNPDE}) is integrated numerically on the rectangle $[-10\pi\,,\,10\pi]^2$ using a finite-difference scheme ($200\times 200$ spatial grid) and 
explicit time-stepping.  Neumann boundary conditions are applied, and initial conditions are chosen so that the spiral tip is away from the edges, so that boundary effects (if any) should be negligible.  

\vspace*{0.2in}
\noindent
{\bf Fully symmetric case}
\vspace*{0.15in}

\begin{figure}[htpb]
\begin{center}
\includegraphics[width=3in]{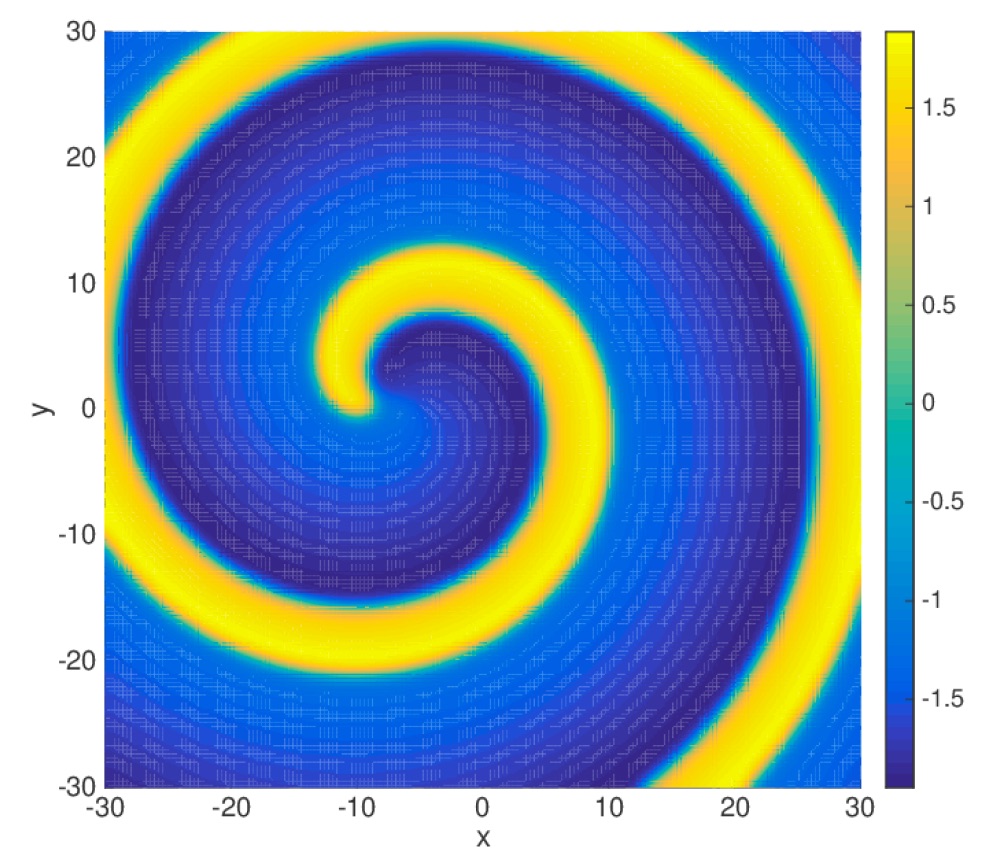}
\includegraphics[width=3in]{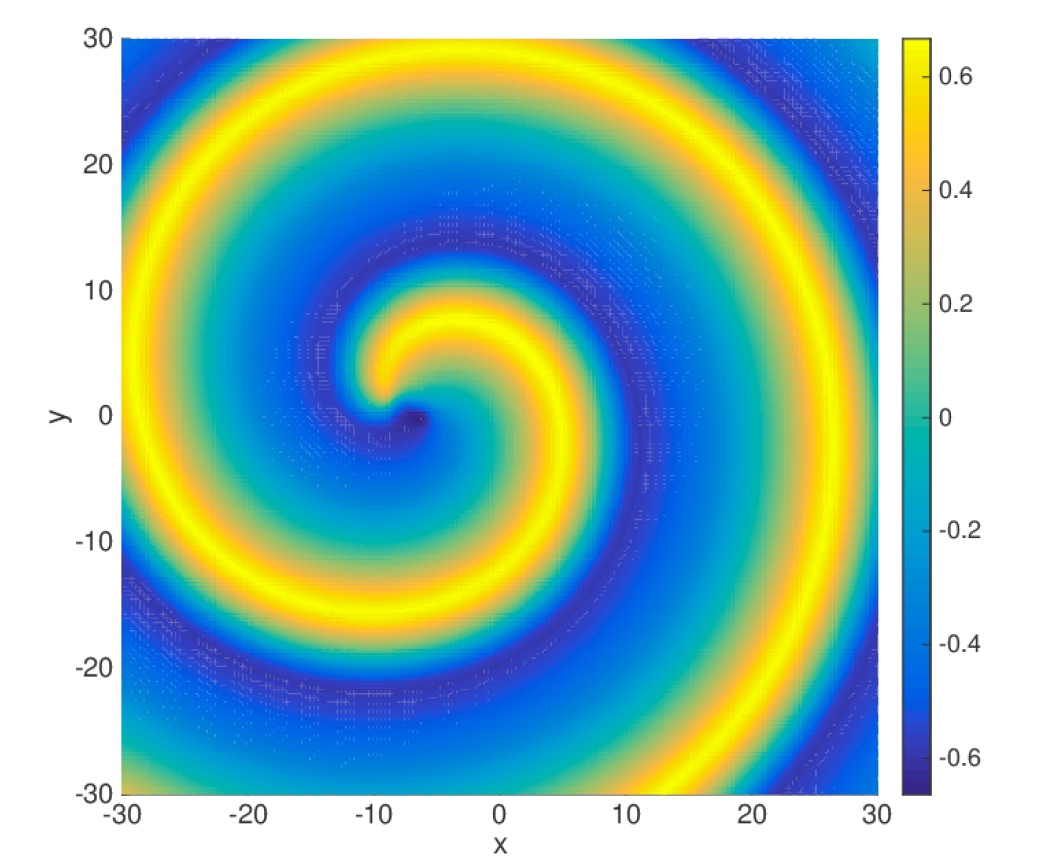}
\end{center}
\caption{Plots of the $u$ (left) and $v$ (right) components at a fixed instant in time for a numerical integration of (\ref{FNPDE}) with $\ep=0$ (homogeneous case), 
$\tau=0.1858$, $\beta=0.755$, $\gamma=0.5$}
\label{fig1}
\end{figure}

\begin{figure}[htpb]
\begin{center}
\includegraphics[width=3.2in]{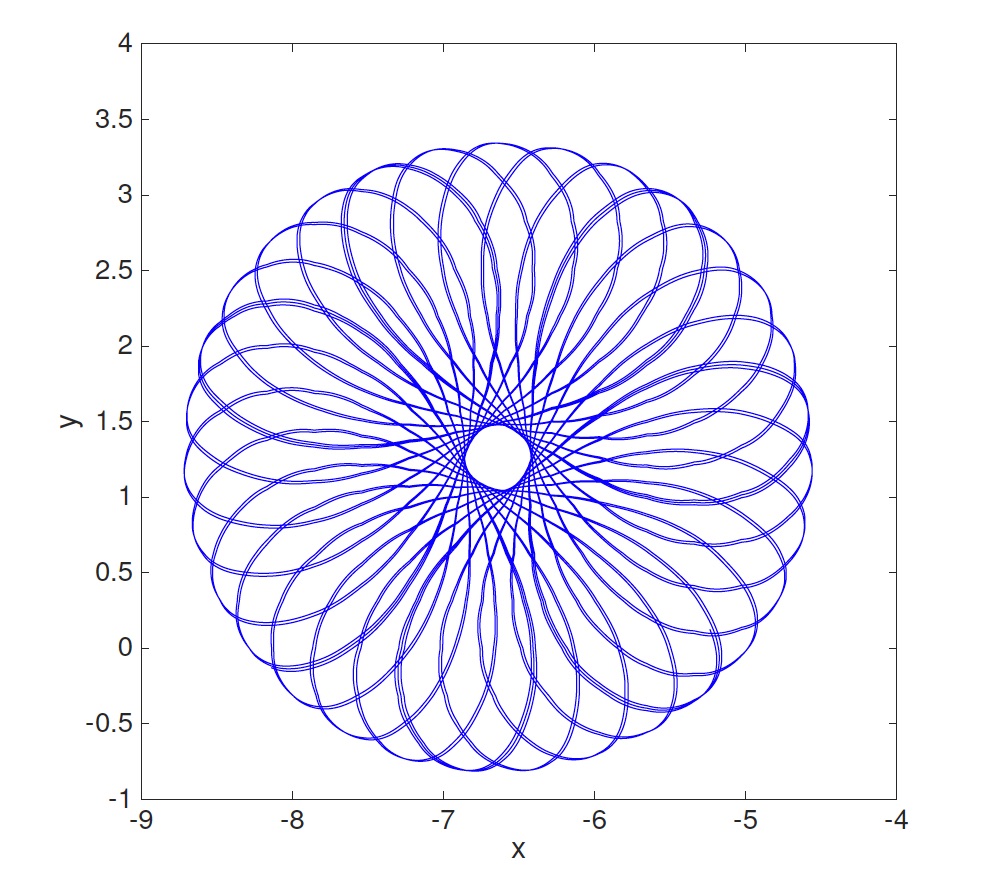}
\includegraphics[width=3.2in]{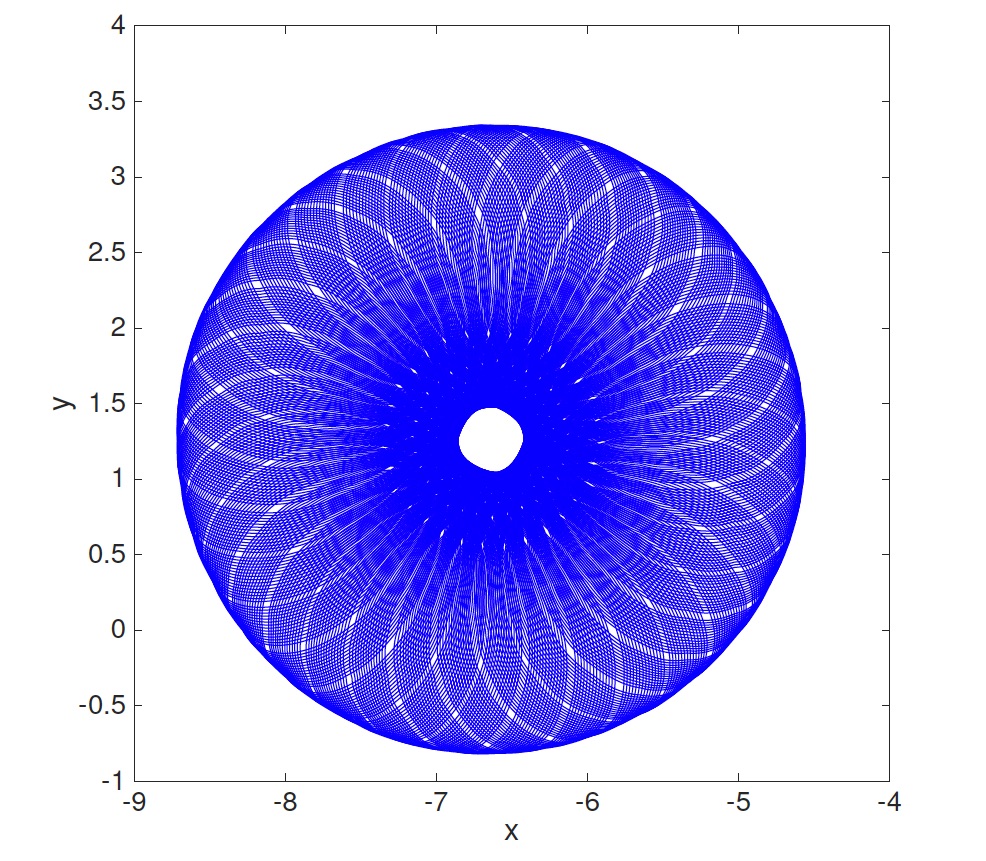}
\end{center}
\caption{Plot in space of the position of the spiral tip for the integration described in Figure \ref{fig1}. The left figure is for a short integration time, and the right figure is on a longer integration time. This is typical of a two-frequency motion, where the frequencies are incommensurate.}
\label{fig2}
\end{figure}

\begin{figure}[htpb]
\begin{center}
\includegraphics[width=3.2in]{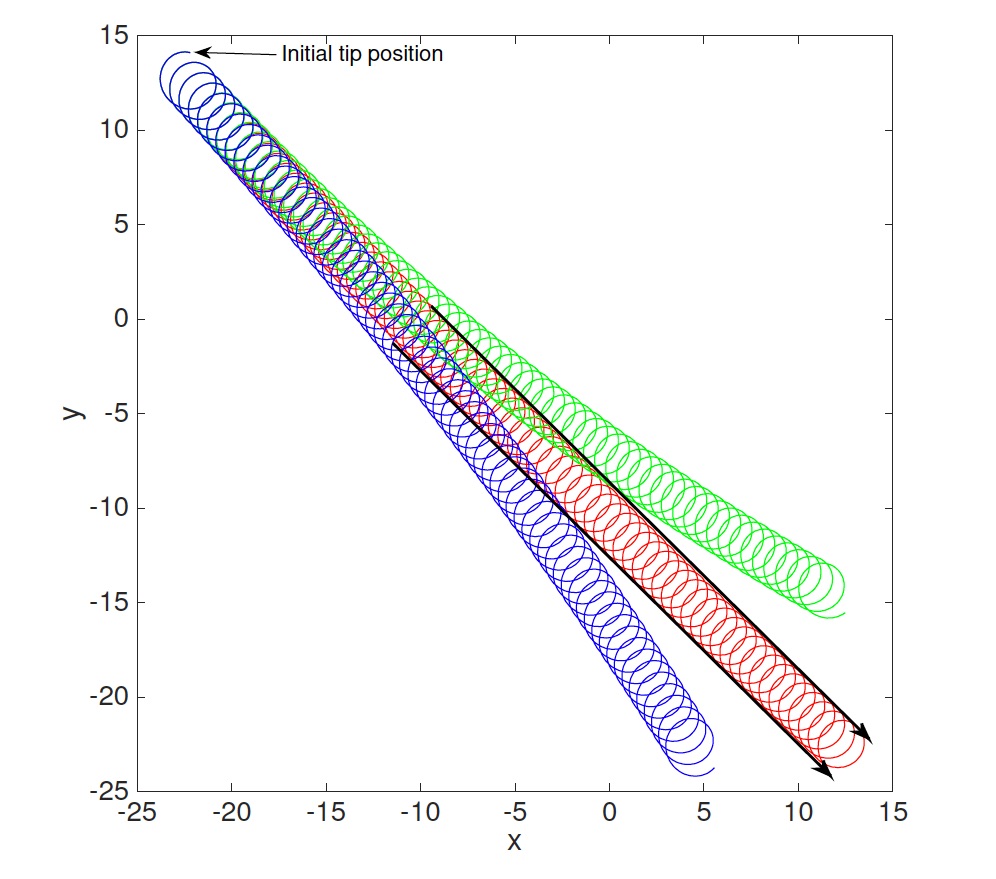}
\end{center}
\caption{Superposition of meander paths for three simulations of (\ref{FNPDE}) with $\ep=0$ and same initial condition in all three cases, which illustrate the transition from meandering with inwards petals (green) to outwards petals (blue) via linear meandering (red).  Parameter values are $\tau=0.26$, $\gamma=0.5$ and $\beta=0.793$ (green), $\beta= 0.79275$ (blue), $\beta=0.792875$ (red).
The parallel black arrows on both sides of the red meander path have been added to the figure to guide the eye.}
\label{fig3}
\end{figure}

\begin{figure}[htpb]
\begin{center}
\includegraphics[width=3.2in]{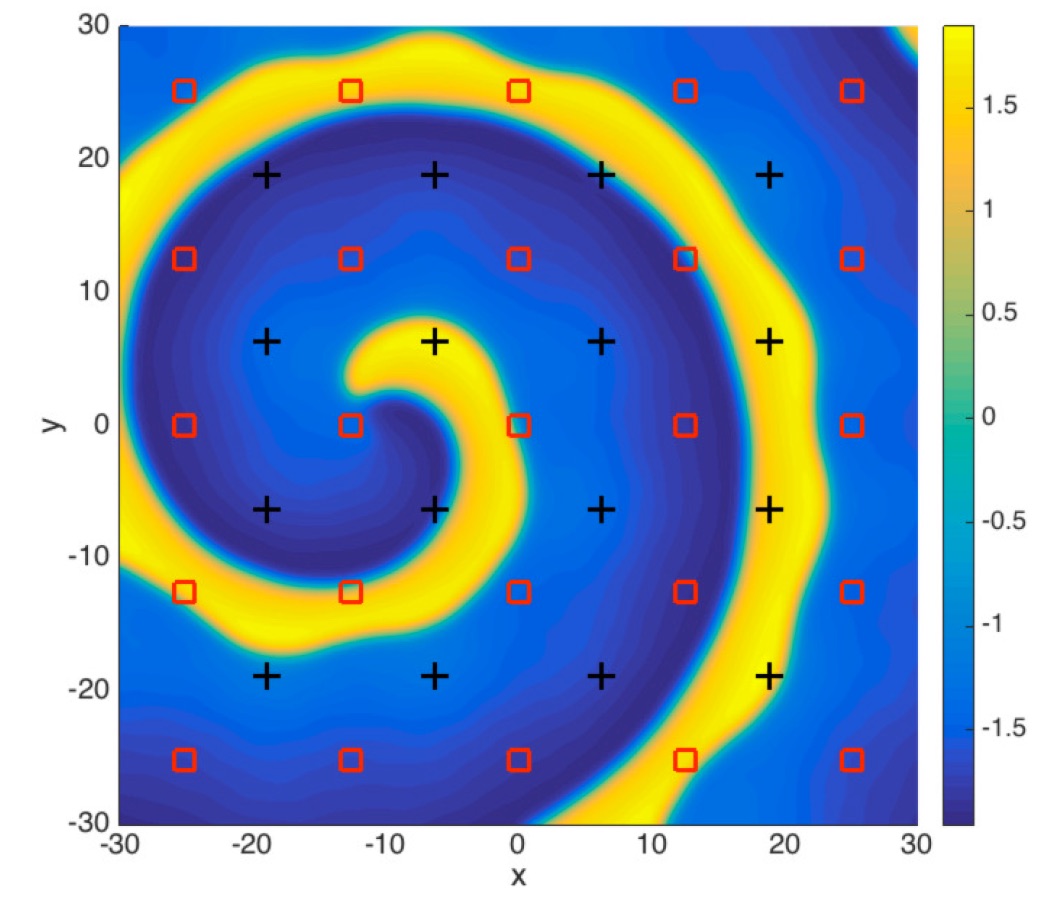}
\end{center}
\caption{Plot of the $u$ component of an integration of (\ref{FNPDE}) with an inhomogeneous perturbation. The effect of the perturbation on the wave profile can be easily observed.  The red squares and black crosses represent the lattice and dual lattice respectively.}
\label{fig4}
\end{figure}

In Figure \ref{fig1}, we give a plot of the $u$ and $v$ components of a typical numerical integration of (\ref{FNPDE}) when $\ep=0$, at a fixed instant in
time, illustrating the spiral shape of the profile.  
To illustrate the meandering dynamics of spiral waves, we will track the position in space of the ``spiral tip'', which we arbitrarily define in this paper as being the intersection of the $u=0$ and
$v=0$ contours.  Note that other definitions of the ``tip'' are possible, and it is well-known \cite{B92,B94,BKT,Roth1,Roth2} that the dynamics are qualitatively similar for different choices of definitions.
In Figure \ref{fig2}, we give the position of the spiral tip for the integration of the homogeneous case described above and illustrated in Figure \ref{fig1}.
We note that the meander path in Figure \ref{fig2} resembles a flower with petals \cite{GLM1} pointing outwards.  Meander paths with inward pointing petals also occur
in (\ref{FNPDE}).  Furthermore, in parameter space, there is a codimension-one surface across which meandering paths with outward petals change continuously into
meandering paths with inward petals via {\em modulated travelling waves}, which are states where the spiral motion is characterized as a composition of a rotation and a linear translation of the center of rotation \cite{BK}.  Figure \ref{fig3} illustrates this phenomenon.

\vspace*{0.2in}
\noindent
{\bf Symmetry-breaking}
\vspace*{0.15in}

For simulations of (\ref{FNPDE}) involving lattice inhomogeneities ($\ep\neq 0$), we typically use as initial condition the $u-v$ state of a homogeneous integration (such as depicted in Figure \ref{fig1}), and
we choose the size of the perturbation terms $|\ep\,g_i|$ small enough so that the wave retains its overall spiral shape.   
For example, in Figure \ref{fig4}, we show the $u$ component
of such an inhomogeneous case.  The effect of the perturbation on the wave profile can be easily observed.  

\subsection{Outline of the paper}

The paper is organized as follows. In section 2, we present the functional analytic framework and hypotheses. The goal is to reduce the problem of studying the effects of 
lattice symmetry-breaking perturbations in (\ref{basicRDPDEpert}) on relative periodic solutions, to that of studying the asymptotic behaviour of solutions of a four-dimensional system of
ordinary differential equations (\ref{cb-sys}) (center bundle equations) possessing symmetry properties related to those of (\ref{basicRDPDEpert}).

One of our most important analytical tools in this paper is Hale's theory of averaging of multiply periodic differential equations and associated
results on the existence of integral manifolds, presented in Theorem 2.3, \S VII.2 of \cite{Hale}.  In section 3, we perform certain changes of
coordinates on (\ref{cb-sys}) in order to transform these equations into forms that are suitable to apply Hale's theorem.  The analysis also depends on the commensurability or
incommensurability of various quantities which will be described later on.  Furthermore, this is related to the classical problem of small divisors \cite{Arnold,Fasano}, which we briefly summarize in section 4.

Our main results are contained in sections 5, 6 and 7 where we study the effects of lattice symmetry-breaking on respectively quasi-periodic meandering waves, on meandering waves whose
meander path is a closed epicycle (we will pay special attention to phase-locking in this case), and on modulated travelling waves.  In each of these 3 chapters, we give the main mathematical
results, give an interpretation for how these mathematical results translate to dynamical features of spiral waves, and then illustrate with numerical results performed on (\ref{FNPDE}) for various choices of inhomogeneity functions $g_i$, and various kinematic parameters $\tau$, $\beta$ and $\gamma$.

We end with a discussion in section 8.  Some of the
more technical proofs are omitted from the main text, and instead presented in Appendix A.

\Section{Preliminaries}

The special Euclidean group, denoted by $\SEtwo$, is the set of all planar translations and rotations.  We parametrize this group as
\begin{equation*}
\SEtwo = \{(\phi, p)\in{\Sone}\times \mathbb{R}^2\}.
\end{equation*}
The action of $\SEtwo$ on a point $z =  \left(\begin{array}{c}x\\y\end{array}\right)\in\mathbb{R}^2,$ is given by
\[
(\phi, p)\cdot z = R_{\phi}\cdot z  + p,
\]
where
\begin{equation}
R_{\phi} = \left(\begin{array}{cr}\cos\phi & -\sin\phi\\ \sin\phi & \cos\phi\end{array}\right), \; p = \left(\begin{array}{c} p_x\\p_y \end{array}\right).
\label{rotmatdef}
\end{equation}
We define $\Sigma$ to be the following subgroup of $\SEtwo$
\begin{equation}
\Sigma = \left\{(\phi, p)\in \SEtwo\,|\, \phi = \frac{n_1\pi}{2}\,(\mbox{\rm mod} \,2\pi), \,p = \left(\begin{array}{c} 2\pi n_2\\2\pi n_3 \end{array}\right), \; n_1, n_2, n_3 \in\mathbb{Z}\right\},
\label{sigdef}
\end{equation}
which represents the symmetry subgroup of a regular square lattice.

Let $X$ be a Banach space, and suppose
\[
a:{\bf SE}(2)\rightarrow GL(X)
\]
is a faithful and isometric representation of ${\bf SE}(2)$ in the space of bounded, invertible linear operators on $X.$ 
The situation we have in mind is that $X$ is a space of functions defined from $\mathbb{R}^2$ into $\mathbb{R}^N$, and the action $a$ on $X$ is
\[
(a(\nu)U)(z) = U(\nu^{-1}\cdot z), \; \nu = (\phi, p)\in \SEtwo.
\]


We consider a semilinear autonomous differential equation on $X$ of the form
\begin{equation}\label{ade}
w_t = {\cal A}w+{\mathscr F}(w) + \ep {\mathscr G}(w,\ep),
\end{equation}
where $\ep \geq 0$ is a small parameter, ${\cal A}$, ${\mathscr F}$ and ${\mathscr G}$ satisfy conditions which guarantee that \eqref{ade} generates a smooth local semiflow $\Phi_{t,\ep} $ on $X$ \cite{Henry},  and ${\mathscr G}$ is bounded.

We assume the following hypothesis on the semiflow $\Phi_{t,\ep}$:
\begin{hypo}
{\rm
\[
\Phi_{t,0}(a(\nu)w) = a(\nu)\Phi_{t,0}(w), \quad \forall w \in X, \, \nu \in \SEtwo, \, t > 0,
\]
and for $\ep > 0,$ we have
\[
\Phi_{t,\ep}(a(\nu)w) = a(\nu)\Phi_{t,\ep}(w), \quad \forall w \in X, \forall t > 0 \Longleftrightarrow \nu\in\Sigma.
\]
}
\end{hypo}
\medskip
This hypothesis means that ${\cal A}+{\mathscr F}$ in \eqref{ade} is $\SEtwo-$equivariant, but ${\mathscr G}$ is only $\Sigma-$equivariant. It is in this sense that we say that ${\mathscr G}$ {\it breaks the $\SEtwo$
symmetry of \eqref{ade}}.


As mentioned in the Introduction, we are interested in how the symmetry-breaking term ${\mathscr G}$ in \eqref{ade} affects the dynamics of meandering waves.  Thus, we will need to assume that
(\ref{ade}) admits such solutions when $\ep=0$.

\begin{hypo}(Existence of normally hyperbolic relative periodic solutions)
{\rm
\newline There exists $w^{\star} \in X$, $T>0$ and $\nu^*\in\SEtwo$ such that
\[
\Phi_{T,0}(w^*)=a(\nu^*)w^*,
\]
and for all $t\in (0,T)$, we have $\Phi_{t,0}(w^*)\notin \{a(\nu)u^*\,|\,\nu\in\SEtwo\}$.  Furthermore, we assume that the set 
$\{\sigma \in \mathbb{C}\;|\;|\sigma| \geq 1\}$
is a spectral set for the linearization  $a(\nu^*)^{-1}D\Phi_{T,0}(w^{\star})$, with projection $P_{\star}$ such that the generalized eigenspace range($P_{\star}$) is four-dimensional (three corresponding to the symmetry eigenvalues, and one corresponding to the flow direction).
\label{hyp2}
}
\end{hypo}

For simplicity, we will only be interested in one-armed spiral waves, so we assume the isotropy subgroup of $w^{\star}$ in Hypothesis \ref{hyp2} is trivial, i.e.\[a(\phi, p)\cdot w^{\star} = w^{\star} \Longleftrightarrow (\phi, p) = (0,0).\]
\begin{rmk}
{\rm
The group element $a(\nu^*)$ in Hypothesis \ref{hyp2} is either a rotation or a translation.  In the former case, the spiral wave evolves in time in such a way that the tip of the
spiral traces out an epicyclic path (with petals either inwards or outwards) - see for example \cite{B92,B94,BKT,Winfree1,Winfree2} who coined such states as {\em meandering}. 
Such a state is illustrated in Figure \ref{fig2}.
In the case where $a(\nu^*)$ is a translation, the path of the spiral tip is a combination of a rotation superimposed with a translation - the center of rotation travels linearly at a constant speed.  These states were described in \cite{BK} as occurring arbitrarily close to a point of resonant Hopf bifurcation from a purely rotating spiral waves. See the red meander path of Figure 
\ref{fig3}
for
an illustration of such a state. 
}
\end{rmk}

Assuming all other hypotheses of the center manifold theorem of \cite{SSW1,SSW2} are satisfied, then for $\ep$ small enough and after a rescaling of time along the orbits of the semiflow, the dynamics of \eqref{ade} near the relative periodic
solution reduce to the following ODE system on the center bundle ${\cal X}$ which is diffeomorphic to the 4-torus $\mathbb{T}^4$ (compare with \cite{CL}):
\begin{equation}\label{cb-sys}
\begin{array}{lll}
\dot{\Psi}& = & R_{\varphi}\left[h_1(\theta) + \varepsilon F_1(\Psi,\varphi,\theta,\varepsilon)\right]\\
\dot{\varphi}& = & \omega + h_2(\theta) + \varepsilon F_2(\Psi,\varphi,\theta,\varepsilon)\\
\dot{\theta}& = & 1,
\end{array}
\end{equation}
where
\[
\left(\begin{array}{c}\Psi\\\varphi\\\theta\end{array}\right)=\left(\begin{array}{c}\psi_1\\\psi_2\\\varphi\\\theta\end{array}\right)\in {\cal X}= \mathbb{T}^4,
\]
$h_{1,2}$ are smooth and $2\pi$-periodic with
\begin{equation}\label{zm}
\langle h_2\rangle\equiv\frac{1}{2\pi}\int_0^{2\pi}h_2(\theta)d\theta = 0,
\end{equation}
$\omega\geq 0$ is a constant real number, $\ep\geq 0$ is a small parameter, and we remind the reader that $R_{\varphi}$ is the rotation matrix (\ref{rotmatdef}).  The variables $\Psi=(\psi_1,\psi_2)$ are the usual translation variables of the center
manifold reduction procedure \cite{CL,GLM1,SSW1,SSW2}, but viewed modulo the lattice of the perturbation, and therefore are periodic with respect to this lattice. It is for this reason that the equations
(\ref{cb-sys}) are viewed as being defined over $\mathbb{T}^4$ instead of the usual $\mathbb{C}\times\mathbb{T}^2$ \cite{LW,LeBlanc}.  When one of the variables $\psi_{1,2}$ advances through
a complete period of $2\pi$, this should be interpreted in physical space as the spiral has advanced to a neighbouring fundamental domain in the spatial grid.

Also in (\ref{cb-sys}), the functions $F_{1,2}$ are assumed to be smooth enough for our purposes, are $2\pi$-periodic in $\psi_{1}$, $\psi_2$, $\varphi$ and $\theta$ and satisfy the
lattice symmetry property
\begin{equation}
\label{symproplat}
F_{1,2}\left(-J\Psi,\varphi+\frac{\pi}{2},\theta,\varepsilon\right)=F_{1,2}(\Psi,\varphi,\theta,\varepsilon),\,\,\,\,\forall\,\left(\begin{array}{c}\Psi\\\varphi\\\theta\end{array}\right)\in\mathbb{T}^4,\,\,0\leq\ep\ll 1,
\end{equation}
where
\[
J\equiv R_{-\frac{\pi}{2}}=\left(\begin{array}{rc}0&1\\-1&0\end{array}\right).
\]

The group $\SEtwo$ acts on the phase space ${\cal X}=\mathbb{T}^4$ of (\ref{cb-sys}) as follows:
\[
(\phi,p)\cdot \left(\begin{array}{c}\Psi\\\varphi\\\theta\end{array}\right)=\left(\begin{array}{c}(R_{\phi}\Psi+p)\,(\mbox{\rm mod}\,[0,2\pi]\times [0,2\pi])\\\varphi+\phi\,\,(\mbox{\rm mod}\,2\pi)\\\theta\end{array}\right).
\]
When $\ep=0$, (\ref{cb-sys}) is equivariant with respect to this action.  However, when $\ep>0$, property (\ref{symproplat}) implies that the symmetry group of (\ref{cb-sys}) is restricted to the lattice subgroup $\Sigma$ defined in (\ref{sigdef}).  Thus, (\ref{cb-sys}) reflects the forced symmetry-breaking property of (\ref{ade}).  When $\ep>0$, if 
$(\Psi(t),\varphi(t),\theta(t))$ is a solution of (\ref{cb-sys}), then so is ${\displaystyle \left(J\Psi(t),\varphi(t)-\frac{\pi}{2},\theta(t)\right)}$.

We end this section with some remarks.
\begin{rmk}

{\rm
\begin{enumerate}
\item[(a)]
The functions $F_{1,2}$, $h_{1,2}$ and the number $\omega$ in (\ref{cb-sys}) obviously depend on the operator ${\cal A}$ and the functions ${\mathscr F}$ and ${\mathscr G}$ in (\ref{ade}), via the center manifold reduction process. 
By studying, as we do in this paper, generic properties of the class of ODEs of the general form (\ref{cb-sys}) without further specifications, we are in a sense studying the 
range of 
possible generic local dynamics
near relative periodic solutions of the infinite-dimensional system (\ref{ade}).  However, from a practical consideration, one would like to be able to obtain a formula that maps between
the space of functions ${\cal A}$, ${\mathscr F}$ and ${\mathscr G}$ in (\ref{ade}), and the space of functions $F_{1,2}$, $h_{1,2}$ and the number $\omega$ in (\ref{cb-sys}).  This is a very difficult problem in general, and
requires information which is usually not readily available (for example, an explicit algebraic representation of the unperturbed meandering spiral wave solution).  With
this limitation in mind, our
goals in this paper will be more modest. We will characterize certain model-independent features of solutions to system (\ref{cb-sys}), and interpret these results in
terms of how a weak lattice symmetry-breaking perturbation generically affects the dynamics of relative periodic solutions of (\ref{ade}).
\item[(b)] Because there has been a rescaling of time to set the coefficient of the $\dot{\theta}$ equation in (\ref{cb-sys}) to 1, the coefficient $\omega$ in (\ref{cb-sys}) in fact represents a
ratio of the two frequencies of the epicyclic meandering in \eqref{ade} when $\ep=0$.
\end{enumerate}
}
\label{remark24}
\end{rmk}

\Section{Preparing the equations for averaging}

In the previous section, we have seen that the local dynamics of (\ref{ade}) near the relative periodic solution reduces to the center bundle ODEs (\ref{cb-sys}) defined on
the 4-torus $\mathbb{T}^4$, where the function $h_2(\theta)$ satisfies the zero-mean condition (\ref{zm}).

Later in the paper, we will want to apply results from 
the theory of averaging and integral manifolds, as set forth in \cite{Hale}, to find invariant periodic solutions and invariant tori to (\ref{cb-sys}) for small $\ep>0$. However, we
must first transform (\ref{cb-sys}) into a standard form suitable to apply Theorem 2.3, \S VII.2 of \cite{Hale} - see also section 3 of \cite{CL}.  As will become evident in the sequel,
the transformations and the resulting transformed system will depend on whether or not $\omega$ is an integer in (\ref{cb-sys}).

\begin{prop}
Suppose $\omega$ in (\ref{cb-sys}) is not an integer. Then under a suitable change of variables, the system (\ref{cb-sys}) is equivalent to
\begin{equation}
\begin{array}{rcl}
\dot{\Psi}&=&\ep R_{\varphi}G_1(\Psi,\varphi,\theta,\ep)\\
\dot{\varphi}&=&\omega+\ep G_2(\Psi,\varphi,\theta,\ep)\\
\dot{\theta}&=&1,
\end{array}
\label{preav_nonint}
\end{equation}
where $G_{1,2}$ are smooth, $2\pi$-periodic in $\psi_1$, $\psi_2$, $\varphi$ and $\theta$, and satisfy the symmetry properties (\ref{symproplat}).
\label{propnonint}
\end{prop}

\begin{prop}
Suppose $\omega$ in (\ref{cb-sys}) is an integer.  Then under a suitable change of variables, the system (\ref{cb-sys}) is equivalent to
\begin{equation}
\begin{array}{rcl}
\dot{\Psi}&=&R_{\varphi}V+\ep R_{\varphi}H_1(\Psi,\varphi,\theta,\ep)\\
\dot{\varphi}&=&\ep H_2(\Psi,\varphi,\theta,\ep)\\
\dot{\theta}&=&1,
\end{array}
\label{preav_int}
\end{equation}
where $V\in\mathbb{R}^2$ is a constant (which is generically non-zero), $H_{1,2}$ are smooth, $2\pi$-periodic in $\psi_1$, $\psi_2$, $\varphi$ and $\theta$, and satisfy the symmetry properties (\ref{symproplat}).
\label{propint}
\end{prop}

The proofs of these propositions are found in Appendix A.

As mentioned above, in the rest of the paper, we will analyze systems (\ref{preav_nonint}) and (\ref{preav_int}) using Theorem 2.3, \S VII.2 of \cite{Hale}.  Throughout, we will assume
that the continuity and smoothness conditions which are required by this theorem are satisfied by ${\cal A}$, ${\mathscr F}$ and ${\mathscr G}$ in (\ref{ade}), and consequently by the vector fields
(\ref{preav_nonint}) and (\ref{preav_int}).  

It is well-known \cite{Arnold,Fasano} that an important technical obstacle to applying ideas from the theory of averaging to a system such as (\ref{preav_nonint}) or (\ref{preav_int}) is the
problem of small divisors. In the next section, we will briefly address this issue.

\Section{Small divisors}

The problem of {\em small divisors} is a classical one in the analysis of multiply periodic dynamical systems, and is related to solving a linear partial differential equation of the form
\[
\omega_1\frac{\partial r}{\partial \phi_1} + \cdots + \omega_{n}\frac{\partial r}{\partial \phi_{n}}=s(\phi_1,\ldots,\phi_{n})
\]
by way of a multiple Fourier series 
\[
r(\phi_1,\ldots,\phi_{n})=\sum_{(m_1,\ldots,m_{n})\in\mathbb{Z}^{n}}\,A_{m_1,\ldots,m_{n}}\,e^{i(m_1\phi_1+\cdots + m_{n}\phi_{n})}
\]
 given
$s(\phi_1,\ldots,\phi_{n})=\sum_{(m_1,\ldots,m_{n})\in\mathbb{Z}^{n}}\,B_{m_1,\ldots,m_{n}}\,e^{i(m_1\phi_1+\cdots + m_{n}\phi_{n})}$ with $B_{0,\ldots,0}=0$.
Formal integration of the PDE leads to 
\[
A_{m_1,\ldots,m_{n}}=\frac{B_{m_{1},\ldots,m_{n}}}{i(m_1\omega_1+\cdots + m_{n}\omega_{n})},\,\,(m_{1},\ldots,m_{n})\in\mathbb{Z}^{n}\setminus \{(0,\ldots,0)\},
\]
which can get uncontrollably large even if $\omega_1,\ldots,\omega_{n}$ are linearly independent over the rationals.  So we can not guarantee the convergence of the series for $r$,
(unless $s$ is a trigonometric polynomial, in which case convergence is not an issue).
We must thus impose restrictions on $\omega_1,\ldots,\omega_{n}$ in the general case.
\begin{mydef}
A vector $\boldsymbol{\Omega}=(\omega_1,\ldots,\omega_{n})\in\mathbb{R}^{n}$, $n\in\mathbb{N}$ is said to satisfy a {\em diophantine condition of constant $\rho>0$ and exponent $\mu\geq n-1$}, and
we write $\boldsymbol{\Omega}\in C_{\rho,\mu}$, if for every $\mathbf{m}=(m_1,\ldots,m_{n})\in\mathbb{Z}^{n}$, we have
\[
| \langle\mathbf{m},\boldsymbol{\Omega}\rangle | \equiv\left|\sum_{j=1}^{n} m_j\omega_j\right|\geq \rho |\mathbf{m}|^{-\mu},
\]
where $|\mathbf{m}|=\sum_{j=1}^{n}|m_j|$.
\label{DiophantineDef}
\end{mydef}
One can show \cite{Fasano} that the Lebesgue measure of $C_{\rho,\mu}$ is positive.

The following is an adaptation of Theorem 12.9 of \cite{Fasano}, where we also use well-known results on the rate of convergence of Fourier coefficients and the degree of smoothness of the corresponding Fourier series \cite{Pereyra}.
\begin{prop}
Suppose $\boldsymbol{\Omega}=(\omega_1,\ldots,\omega_{n})\in\mathbb{R}^{n}$ is such that its components are linearly independent over the rationals, and consider the partial differential equation on the torus $\mathbb{T}^{n}$
\begin{equation}
\sum_{k=1}^{n}\,\omega_{k}\,\frac{\partial r}{\partial \phi_{k}}(\mathbf{x},\phi_1,\ldots,\phi_{n}) = s(\mathbf{x},\phi_1,\ldots,\phi_{n})
\label{genPDE}
\end{equation}
where $\mathbf{x}\in K\subset\mathbb{R}^{q}$, $K$ compact.
\begin{enumerate}
\item[(a)] If $s(\mathbf{x},\boldsymbol{\phi})=s(\mathbf{x},\phi_1,\ldots,\phi_{n})$ is a trigonometric polynomial of the form
\[
s(\mathbf{x},\boldsymbol{\phi})=\sum_{m_1=-N_1}^{N_1}\cdots\sum_{m_{n}=-N_1}^{N_1}\,B_{\mathbf{m}}(\mathbf{x})e^{i\langle\mathbf{m},\boldsymbol{\phi}\rangle}
\]
with $B_{\mathbf{0}}=0$, then the trigonometric polynomial
\[
r(\mathbf{x},\boldsymbol{\phi})=\sum_{m_1=-N_1}^{N_1}\cdots\sum_{m_{n}=-N_1}^{N_1}\,A_{\mathbf{m}}(\mathbf{x})e^{i\langle\mathbf{m},\boldsymbol{\phi}\rangle},
\]
where
$A_{\mathbf{0}}=0$ and ${\displaystyle A_{\mathbf{m}}(\mathbf{x})=\frac{B_{\mathbf{m}}(\mathbf{x})}{i\langle\mathbf{m},\boldsymbol{\Omega}\rangle}},\,\,\mathbf{m}\neq \mathbf{0}$, is a solution
to (\ref{genPDE}).
\item[(b)] If ${\displaystyle s(\mathbf{x},\boldsymbol{\phi})=\sum_{\mathbf{m}\in\mathbb{Z}^{n}}\,B_{\mathbf{m}}(\mathbf{x})e^{i\langle\mathbf{m},\boldsymbol{\phi}\rangle}}$ (with $B_{\mathbf{0}} =\mathbf{0}$) is smooth, but not a trigonometric polynomial, and if $\boldsymbol{\Omega}\in C_{\rho,\mu}$ (for some $\rho>0$ and $\mu\geq n-1$), then the function
\[
r(\mathbf{x},\boldsymbol{\phi})=\sum_{\mathbf{m}\in\mathbb{Z}^{n}\setminus\{\mathbf{0}\}}\,\frac{B_{\mathbf{m}}(\mathbf{x})}{i\langle\mathbf{m},\boldsymbol{\Omega}\rangle}e^{i\langle\mathbf{m},\boldsymbol{\phi}\rangle}
\]
is smooth (with degree of smoothness determined by $\mu$), and is a solution to (\ref{genPDE}).
\end{enumerate}
\label{DiophantineProp}
\end{prop}

\Section{The case where $\omega$ is irrational in (\ref{preav_nonint})}

Performing the near identity change of variables $\Psi=\hat{\Psi}+\varepsilon r(\hat{\Psi},\varphi,\theta)$ transforms (\ref{preav_nonint}) into (upon dropping the hats)
\begin{equation}
\begin{array}{rcl}
\dot{\Psi}&=&\ep {\cal G}(\Psi)+\ep \left[s(\Psi,\varphi,\theta)-\omega\frac{\partial r}{\partial\varphi}(\Psi,\varphi,\theta)-\frac{\partial r}{\partial\theta}(\Psi,\varphi,\theta)\right]+\ep^2 R_{\varphi}G_3(\Psi,\varphi,\theta,\ep)\\
\dot{\varphi}&=&\omega+\ep G_4(\Psi,\varphi,\theta,\ep)\\
\dot{\theta}&=&1,
\end{array}
\label{preav_nonint_1}
\end{equation}
where
\begin{equation}
{\cal G}(\Psi)=\frac{1}{(2\pi)^2}\int_{\mathbb{T}^2}\,R_{\varphi}G_1(\Psi,\varphi,\theta,0)\,d\varphi d\theta,
\label{GDEF}
\end{equation}
and
\begin{equation}
s(\Psi,\varphi,\theta)=R_{\varphi}G_1(\Psi,\varphi,\theta,0)-{\cal G}(\Psi).
\label{newsdef}
\end{equation}
The following is a straightforward consequence of Proposition \ref{DiophantineProp}, 
\begin{prop}
If $s$ in (\ref{preav_nonint_1}) and (\ref{newsdef}) is a trigonometric polynomial, or if $\boldsymbol{\Omega}=(\omega,1)$ satisfies Definition \ref{DiophantineDef}, then
a smooth function $r$ can be chosen so that (\ref{preav_nonint_1}) simplifies to
\begin{equation}
\begin{array}{rcl}
\dot{\Psi}&=&\ep {\cal G}(\Psi)+\ep^2 R_{\varphi}G_3(\Psi,\varphi,\theta,\ep)\\
\dot{\varphi}&=&\omega+\ep G_4(\Psi,\varphi,\theta,\ep)\\
\dot{\theta}&=&1,
\end{array}
\label{preav_nonint_2}
\end{equation}
where ${\cal G}$ is as in (\ref{GDEF}).
Moreover, $G_{3,4}$ satisfy the symmetry properties (\ref{symproplat}), and ${\cal G}(J\Psi)=J{\cal G}(\Psi)$.
\end{prop}

The two-dimensional $\mathbb{Z}_4-$equivariant system
 \begin{equation}\label{eq:av}
\dot{\Psi} = \ep\,{\cal G}(\Psi),
\end{equation}
has an equilibrium at $\Psi = 0$, and any non-trivial equilibria occur as a conjugate set
\[
\left\{J^k\Psi^{\star}, \; k=0, 1, 2, 3\right\}.
\]
If $\left\{\Psi(t) |\; 0\leq t \leq T\right\}$ is a $T-$periodic orbit of \eqref{eq:av}, then either  $\left\{J^k\Psi(t) |\; 0\leq t \leq T\right\}$ is a distinct periodic orbit for $k= 1, 2, 3,$ or $\left\{J\Psi(t) |\; 0\leq t \leq T\right\}$ coincides with the orbit $\left\{\Psi(t) |\; 0\leq t \leq T\right\}$. In this latter case, the periodic solution has one of the following spatial-temporal symmetries
\begin{equation}\label{sp-temp-sym}
\Psi(t - T/4) = \pm J\Psi(t).
\end{equation}

The main result of this section is the following:
\begin{theo}\label{theo3}
Let us consider the equation \eqref{eq:av}.
\begin{enumerate}
\item[(i)] Suppose the equation \eqref{eq:av} has an equilibrium at $\Psi^{\star}$ which is linearly stable (resp. unstable). Then for $\ep > 0$ small enough, the system \eqref{preav_nonint_2}
has linearly stable (resp. unstable) invariant two-tori represented as
\begin{equation}\label{2T-Sol}
\Psi = J^k\left(\Psi^{\star} + \sqrt{\ep}\sigma_{\Psi^{\star}}\left(\varphi + \frac{k\pi}{2},\theta,\ep\right)\right), \; k = 0, 1, 2, 3,
\end{equation}
where the smooth function $\sigma_{\Psi^{\star}}$ is such that $\sigma_{\Psi^{\star}} \rightarrow 0$ as $\ep \rightarrow 0.$
If $\Psi^{\star} = 0$, then the torus has the $\mathbb{Z}_4$-symmetry
\begin{equation}\label{2T-Sym}
\sigma_{0}\left(\varphi - \frac{\pi}{2},\theta,\ep\right) = J\sigma_{0}(\varphi,\theta,\ep).
\end{equation}
\item[(ii)] Suppose  $\left\{\Psi^{\star}(t) |\; 0\leq t \leq T\right\}$ is a non-trivial periodic orbit of \eqref{eq:av} and that this orbit is linear stable (resp. unstable). Then for $\ep > 0$ small enough, the system  \eqref{preav_nonint_2} has linearly stable (resp. unstable) invariant three-tori represented as 
\begin{equation}\label{3T-Sol}
\begin{array}{ll}
\Psi = J^k\left(\Psi^{\star}(\eta + kT/4) \right.&+\left. \sqrt{\ep}\,\Sigma_{\Psi^{\star}}\left(\eta + kT/4, \varphi + {k\pi}/{2},\theta,\ep\right)\right), \, k = 0, 1, 2, 3,\\[0.15in]
&(\eta,\varphi,\theta)\in [0,T]\times [0,2\pi]\times [0,2\pi],
\end{array}
\end{equation}
where the smooth function $\Sigma_{\Psi^{\star}}$ is such that $\Sigma_{\Psi^{\star}} \rightarrow 0$ as $\ep \rightarrow 0.$
Furthermore, if $\Psi^{\star}$ satisfies the spatial-temporal symmetry \eqref{sp-temp-sym},  then the invariant three-torus is such that
\begin{equation}\label{3T-Sym}
\Sigma_{\Psi^{\star}}(\eta - T/4, \varphi \mp {\pi}/{2},\theta,\ep) = \pm J\Sigma_{\Psi^{\star}}(\eta, \varphi,\theta,\ep).
\end{equation}
\end{enumerate}
\end{theo}
\proof
We first prove item (i). Setting 
\begin{equation}
\Psi=\Psi^*+\sqrt{\ep}\hat{\Psi}
\label{cov1}
\end{equation}
transforms (\ref{preav_nonint_2}) into the following (upon dropping the hats)
\begin{equation}
\begin{array}{rcl}
\dot{\Psi}&=&\ep D_{\Psi}{\cal G}(\Psi^*)\Psi+\ep^{\frac{3}{2}} R_{\varphi}{\cal R}_1(\Psi,\varphi,\theta,\ep)\\
\dot{\varphi}&=&\omega+\ep {\cal R}_2(\Psi,\varphi,\theta,\ep)\\
\dot{\theta}&=&1,
\end{array}
\label{preav_nonint_4}
\end{equation}
where ${\cal R}_{1,2}$ are smooth and $2\pi$-periodic in $\varphi$ and in $\theta$.  System (\ref{preav_nonint_4}) is in the appropriate form to apply 
Theorem 2.3, \S VII.2 of \cite{Hale}, and
conclude that (\ref{preav_nonint_4}) has an invariant two-torus of the form $\Psi=\sigma_{\Psi^*}(\varphi,\theta,\varepsilon)$ for all $\ep>0$ sufficiently small, with
$\sigma_{\Psi^*}\rightarrow 0$ as $\ep\rightarrow 0$.
Using (\ref{cov1}), we get that (\ref{2T-Sol}) with $k=0$ is an invariant two-torus for (\ref{preav_nonint_2}).

For the conjugate equilibrium $J\Psi^*$ of (\ref{eq:av}), we replace (\ref{cov1}) by
\[
\Psi=J(\Psi^*+\sqrt{\ep}\hat{\Psi}).
\]
Setting $\varphi=\hat{\varphi}-\pi/2$ and using the equivariance properties of ${\cal R}_1$ and ${\cal R}_2$, (\ref{preav_nonint_2}) also transforms into (\ref{preav_nonint_4}) upon
dropping the hats.  Therefore, $\Psi=J\Psi^*+J\sqrt{\ep}\sigma_{\Psi^*}(\varphi+\pi/2,\theta,\ep)$ is also an invariant two-torus for (\ref{preav_nonint_2}), distinct from the
previous one if $\Psi^*\neq 0$, but if $\Psi^*=0$, we get (\ref{2T-Sym}). This ends the proof of item (i).

To prove (ii), we introduce a local coordinate system $(b,\eta)$ near the periodic orbit (see \cite{Hale} and \cite{CL}) defined by
\begin{equation}
\Psi=\Psi^*(\eta)+b\,J\frac{\dot{\Psi}^*(\eta)}{||\dot{\Psi}^*(\eta)||}
\label{local_coords}
\end{equation}
which transforms (\ref{preav_nonint_2}) into
\begin{equation}
\begin{array}{rcl}
\dot{b}&=&\ep A(\eta,b)b+\ep^2 {\cal U}(b,\eta,\varphi,\theta,\ep)\\
\dot{\eta}&=&\ep (1+B(\eta,b)b)+\ep^2 {\cal V}(b,\eta,\varphi,\theta,\ep)\\
\dot{\varphi}&=&\omega+\ep {\cal W}(b,\eta,\varphi,\theta,\ep)\\
\dot{\theta}&=&1,
\end{array}
\label{3Tavs}
\end{equation}
where ${\cal U}$, ${\cal V}$ and ${\cal W}$ are smooth, $T$-periodic in $\eta$, and $2\pi$-periodic in $\varphi$ and $\theta$, and where
\[
A(\eta,0)=\mbox{\rm trace}(D_{\Psi}{\cal G}(\Psi^*(\eta)))-\frac{d}{d\eta}\ln (||\dot{\Psi}^*(\eta)||).
\] 
If we define
\[
\delta=\frac{1}{T}\int_0^T\,A(s,0)\,ds=\frac{1}{T}\int_0^T\,\mbox{\rm trace}(D_{\Psi}{\cal G}(\Psi^*(s)))\,ds,
\]
then $\delta<0$ (resp., $\delta>0$) if the limit cycle $\Psi^*(t)$ is linearly stable (resp., unstable). Performing the periodic change of variable
\[
b=\sqrt{\ep}\,\xi e^{\int_0^{\eta}\,(A(s,0)-\delta)\,ds}
\]
transforms (\ref{3Tavs}) into
\begin{equation}
\begin{array}{rcl}
\dot{\xi}&=&\ep \delta \xi+\ep^{\frac{3}{2}} {U}(\xi,\eta,\varphi,\theta,\ep)\\
\dot{\eta}&=&\ep +\ep^{\frac{3}{2}} {V}(\xi,\eta,\varphi,\theta,\ep)\\
\dot{\varphi}&=&\omega+\ep {W}(\xi,\eta,\varphi,\theta,\ep)\\
\dot{\theta}&=&1,
\end{array}
\label{3Tavs2}
\end{equation}
which is in the proper form to apply Theorem 2.3, \S VII.2 of \cite{Hale}, and conclude that (\ref{3Tavs2}) admits an invariant three-torus
\[
\xi=\sigma_{\Psi^*}(\eta,\varphi,\theta,\ep),\,\,\,\,\,\sigma_{\Psi^*}\rightarrow 0,\,\,\mbox{\rm as}\,\,\ep\rightarrow 0.
\]
Using (\ref{local_coords}), we conclude that (\ref{preav_nonint_2}) admits an invariant three-torus
\[
\Psi=\Psi^*(\eta)+\sqrt{\ep}\Sigma_{\Psi^*}(\eta,\varphi,\theta,\ep),
\]
where
\[
\Sigma_{\Psi*}(\eta,\varphi,\theta,\ep)=\sigma_{\Psi^*}(\eta,\varphi,\theta,\ep)e^{\int_0^{\eta}\,(A(s,0)-\delta)\,ds}\frac{J\dot{\Psi}^*(\eta)}{||\dot{\Psi}^*(\eta)||}.
\]
We have thus established (\ref{3T-Sol}) in the case $k=0$.  The other cases for $k$ in (\ref{3T-Sol}) and the spatio-temporal symmetry property (\ref{3T-Sym}) follow using similar
arguments to those used in the proof of item (i) above.
\hfill\qed

\noindent
{\bf Interpretation}:  Two-frequency meandering spiral waves have been observed in RDPDE systems such as (\ref{FNPDE}) when $\ep=0$ \cite{BK}.  If the two-frequencies are incommensurate, then the meander path typically resembles that in Figure \ref{fig2}, and this path densely fills a circular annulus.
Theorem \ref{theo3} characterizes the simplest ways in which these meandering waves typically behave under a generic lattice symmetry-breaking perturbation, assuming the frequency ratio satisfies a Diophantine condition.

One possibility is that the two-frequency meandering wave drifts and anchors at a point in space that may or may not be a lattice point.  If the point of anchoring is a point of the lattice, then the meander path will typically lose the full circular symmetry of Figure \ref{fig2}, but retains a four-fold rotational symmetry.  If the point of anchoring is not related to the lattice, then we expect that
the meander path doesn't retain any rotational symmetries.  Furthermore, to any such anchored meander paths, there are conjugate anchored meander paths related to the first one by the symmetries of the lattice.

Another possibility is that the lattice symmetry-breaking induces a third (slow) frequency in the meandering motion of the spiral.  Again, the resultant three-frequency motion will be anchored at some
point in space.  If this anchoring point is a lattice point, then the three-frequency meander path has overall four-fold rotational symmetry, otherwise it retains no rotational symmetries.
Again, lattice symmetries lead
to conjugate three-frequency meandering states.

\begin{figure}[htpb]
\begin{center}
\includegraphics[width=3.2in]{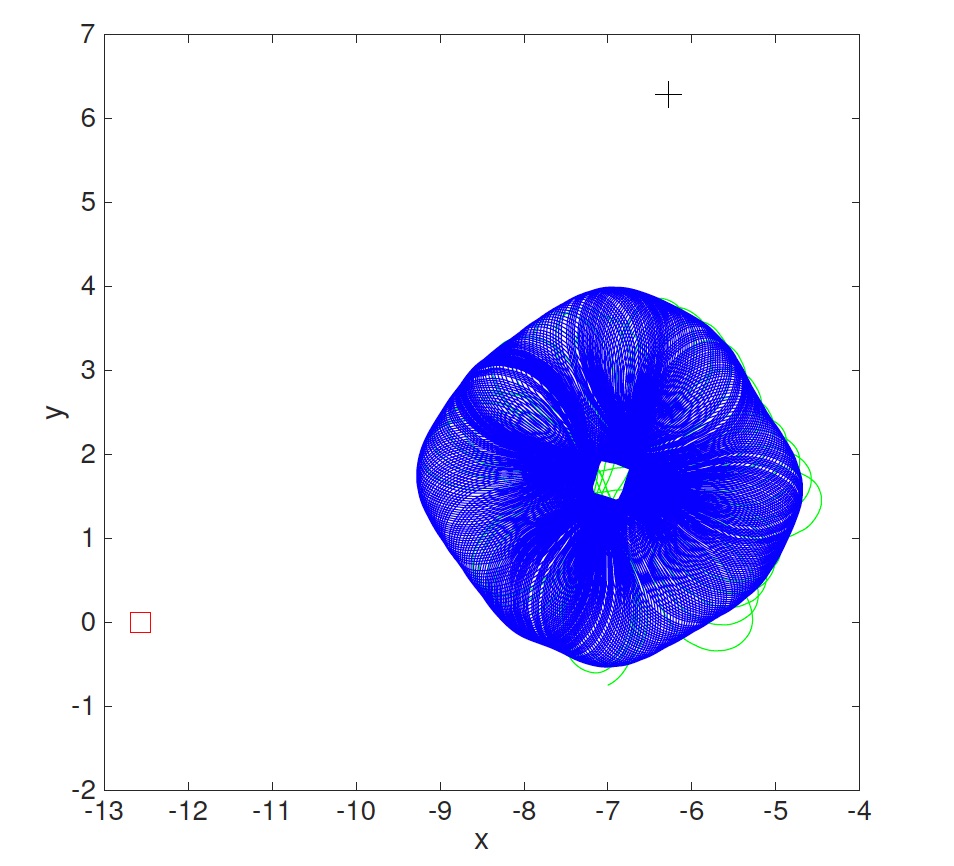}
\end{center}
\caption{A two-frequency anchored meandering path for a spiral wave in (\ref{FNPDE}) with inhomogeneity data (\ref{coeffsinhirrat1}) and kinetic parameters 
$\tau=0.1858$, $\beta=0.755$, $\gamma=0.5$.  The red square and black cross are there for reference purposes to illustrate points of the lattice and the dual lattice respectively. The transient
is in green, and the final anchored state is in blue.
}
\label{irratfig1}
\end{figure}

\begin{figure}[htpb]
\begin{center}
\includegraphics[width=3.2in]{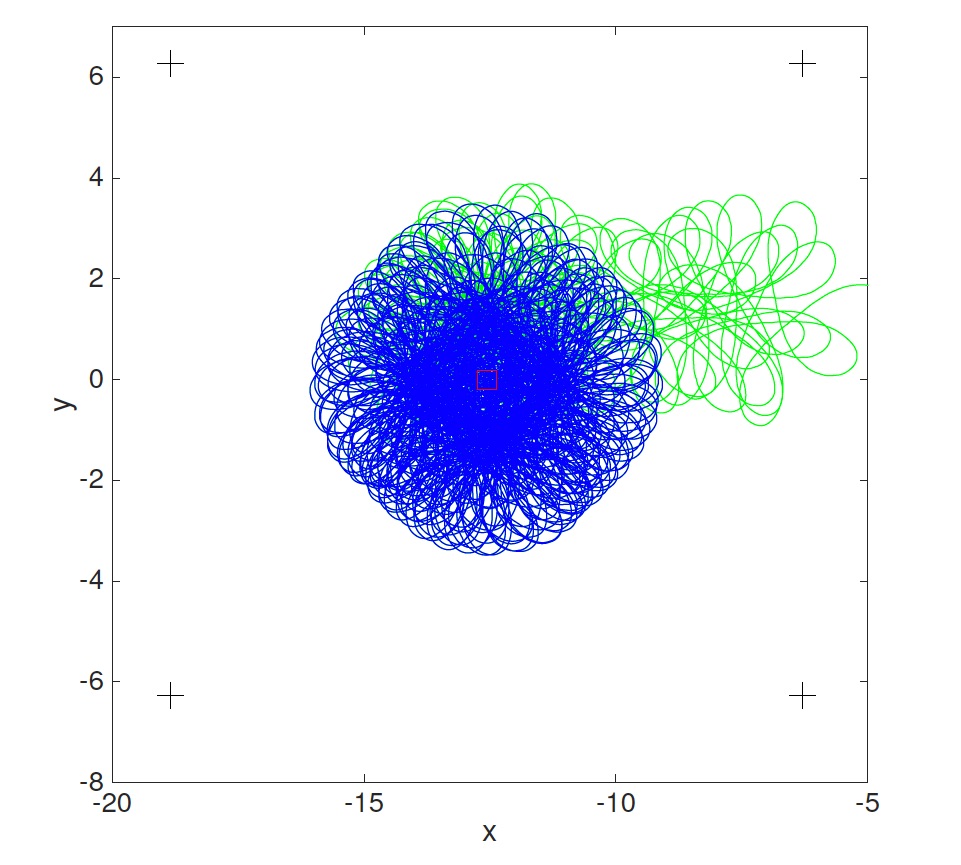}
\end{center}
\caption{A three-frequency anchored meandering path for a spiral wave in (\ref{FNPDE}) with inhomogeneity data (\ref{coeffsinhirrat2}) and kinetic parameters 
$\tau=0.1858$, $\beta=0.755$, $\gamma=0.5$.  
The red square and black crosses are there for reference purposes to illustrate points of the lattice and the dual lattice respectively. The transient
is in green, and the final anchored state is in blue.
}
\label{irratfig2}
\end{figure}

In Figures \ref{irratfig1} and \ref{irratfig2}, we show the meander path of a spiral wave solution observed in a numerical simulations of (\ref{FNPDE}) with inhomogeneity data (\ref{gspecif}) given respectively by
\begin{equation}
\ep=0.01,\,\,A_1=-0.7,\,\,A_2=0.14,\,\,B_1=-2.5,\,\,B_2=-0.5,\,\,C_1=0.5,\,\,C_2=1.5\,\,\mbox{\rm for Figure \ref{irratfig1}}
\label{coeffsinhirrat1}
\end{equation}
and
\begin{equation}
\ep=0.01,\,\,A_1=-0.9,\,\,A_2=0.46,\,\,B_1=-2.5,\,\,B_2=0.5,\,\,C_1=-0.5,\,\,C_2=1.5\,\,\mbox{\rm for Figure \ref{irratfig2}}.
\label{coeffsinhirrat2}
\end{equation}
We have chosen the same kinetic parameters $\tau=0.1858$, $\beta=0.755$, $\gamma=0.5$ as for the simulation of the fully Euclidean case illustrated in Figure \ref{fig2}.  After transients have died out, one observes a two-frequency meandering wave which is anchored at a point that is not a lattice point in Figure \ref{irratfig1}, and a three-frequency meandering wave which is anchored at a lattice point in
Figure \ref{irratfig2}

\Section{Phase-locking}

When $\omega=k/\ell$ is rational (but not an integer) in (\ref{preav_nonint}) (with $\mbox{\rm gcd}(k,\ell)=1$), then as we will see in this section, phase-locking of meandering waves can occur.  We prepare (\ref{preav_nonint}) for averaging by setting
$\hat{\varphi}=\varphi-k\theta/\ell$, which yields (after dropping the hats)
\begin{equation}
\begin{array}{rcl}
\dot{\Psi}&=&\ep R_{\varphi}R_{k\theta/\ell}{G}_1(\Psi,\varphi+k\theta/\ell,\theta,\ep)\\
\dot{\varphi}&=&\ep {G}_2(\Psi,\varphi+k\theta/\ell,\theta,\ep)\\
\dot{\theta}&=&1.
\end{array}
\label{preav_nonint_rational}
\end{equation}
We then define the averages
\begin{equation}
\begin{array}{c}
{\displaystyle{\cal G}_1(\Psi,\varphi)=\frac{1}{2\pi\ell}\int_0^{2\pi\ell}\,R_{\varphi}R_{k\theta/\ell}G_1(\Psi,\varphi+k\theta/\ell,\theta,0)\,d\theta}\\[0.2in]
{\displaystyle{\cal G}_2(\Psi,\varphi)=\frac{1}{2\pi\ell}\int_0^{2\pi\ell}\,G_2(\Psi,\varphi+k\theta/\ell,\theta,0)\,d\theta}.
\end{array}
\label{averages_rat}
\end{equation}
The following is a straightforward consequence of the standard (singly periodic) theory of averaging \cite{GH,Hale}
\begin{theo}
Let $(\Psi_0,\varphi_0)$ be a hyperbolic equilibrium point of the averaged equations
\begin{equation}
\begin{array}{rcl}
\dot{\Psi}&=&\ep {\cal G}_1(\Psi,\varphi)\\
\dot{\varphi}&=&\ep {\cal G}_2(\Psi,\varphi).
\end{array}
\label{aveqsrat}
\end{equation}
\label{mainthm_avg_rational}
Then for all small enough $\varepsilon>0$, the system (\ref{preav_nonint_rational}) has a $2\pi\ell$-periodic solution represented as
\begin{equation}
\Psi=f^{1}_{(\Psi_0,\varphi_0)}(\theta,\ep),\,\,\,\,\,\,\,\varphi=f^{2}_{(\Psi_0,\varphi_0)}(\theta,\ep)
\label{solpredav}
\end{equation}
which tend to $\Psi_0$ and $\varphi_0$ respectively when $\ep\rightarrow 0$.
Furthermore, the stability of the periodic solution is the same as the stability of the equilibrium point $(\Psi_0,\varphi_0)$ in (\ref{aveqsrat}).
\label{thmmainphaselock}
\end{theo}

\noindent
{\bf Interpretation}:
Theorem \ref{thmmainphaselock} implies that if $(\Psi_0,\varphi_0)$ is a linearly stable equilibrium point for (\ref{aveqsrat}), then
(\ref{preav_nonint_rational}) has an asymptotically stable
$2\pi\ell$-periodic solution 
\[
(\Psi(t),\varphi(t))=(f^1_{(\Psi_0,\varphi_0)}(t,\ep),f^2_{(\Psi_0,\varphi_0)}(t,\ep)).
\]
Working backwards through the changes of coordinates that transformed (\ref{cb-sys}) into (\ref{preav_nonint_rational}) (see Appendix A for details), the linearly stable $2\pi\ell$-periodic solution (\ref{solpredav})
 of (\ref{preav_nonint_rational}) corresponds to a $2\pi\ell$-periodic solution of (\ref{cb-sys}) which is such that
\begin{equation}
\Psi(t,\ep)=f^1_{(\Psi_0,\varphi_0)}(t,\ep)+R_{f^2_{(\Psi_0,\varphi_0)}(t,\ep)}R_{kt/\ell}R_{\int_0^t\,h_2(s)\,ds}\,{\cal S}(t),
\label{workback}
\end{equation}
where ${\cal S}(t)$ is $2\pi$-periodic (see Lemma \ref{lemQ} and equation (\ref{fundcoc})).  
The term 
\[
R_{f^2_{(\Psi_0,\varphi_0)}(t,\ep)}R_{kt/\ell}R_{\int_0^t\,h_2(s)\,ds}\,{\cal S}(t)
\]
represents a compound motion of a $2\pi$-periodic term
$R_{\int_0^t\,h_2(s)\,ds}\,{\cal S}(t)$, a $2\pi\ell/k$-periodic term $R_{kt/\ell}$ and a $2\pi\ell$-periodic term $R_{f^2_{(\Psi_0,\varphi_0)}(t,\ep)}$.
Thus, this $2\pi\ell$-periodic solution corresponds to a two-frequency solution for (\ref{cb-sys}), where the two frequencies are in ratio $k/\ell$.
If such a solution represents the dynamics of a spiral wave solution to (\ref{ade}) via the center manifold reduction (\ref{cb-sys}), then
the meander path of the spiral tip is
closed and is composed of two-frequencies which are in ratio $k/\ell$.  Moreover, because of the linear stability of this solution, then it is stable to small perturbations in (\ref{ade}).
This means that contrary to the fully Euclidean case, lattice symmetry-breaking may lead to phase-locking of meandering spiral waves.  

\begin{figure}[htpb]
\begin{center}
\includegraphics[width=3.2in]{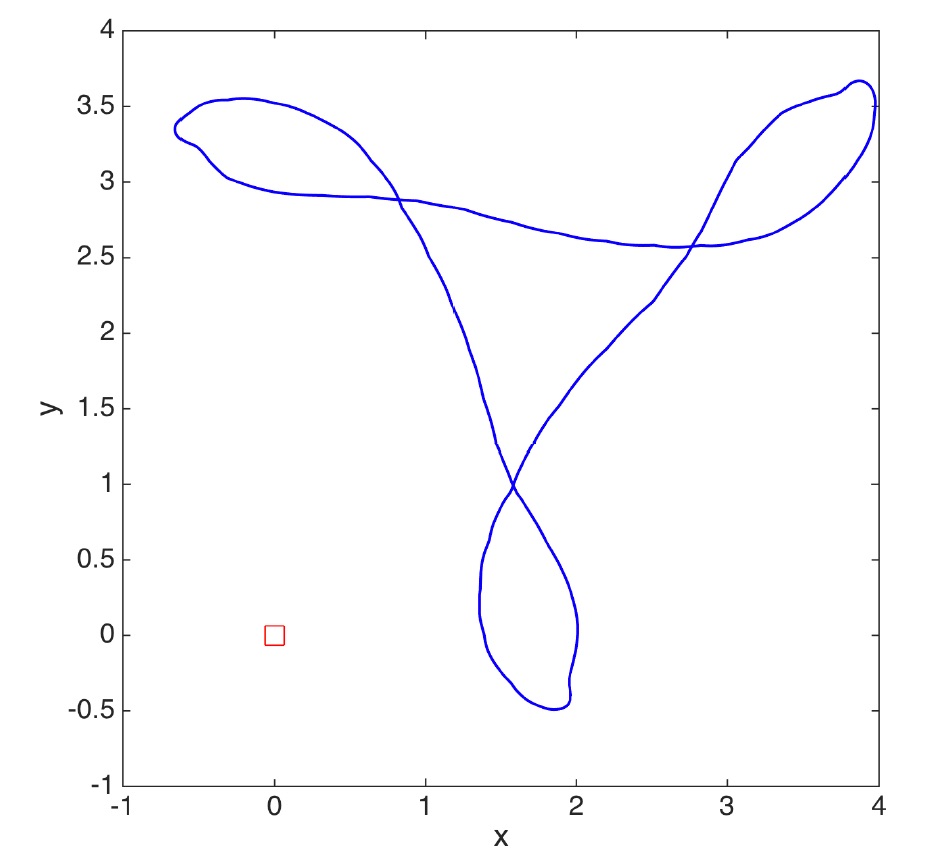}
\includegraphics[width=3.2in]{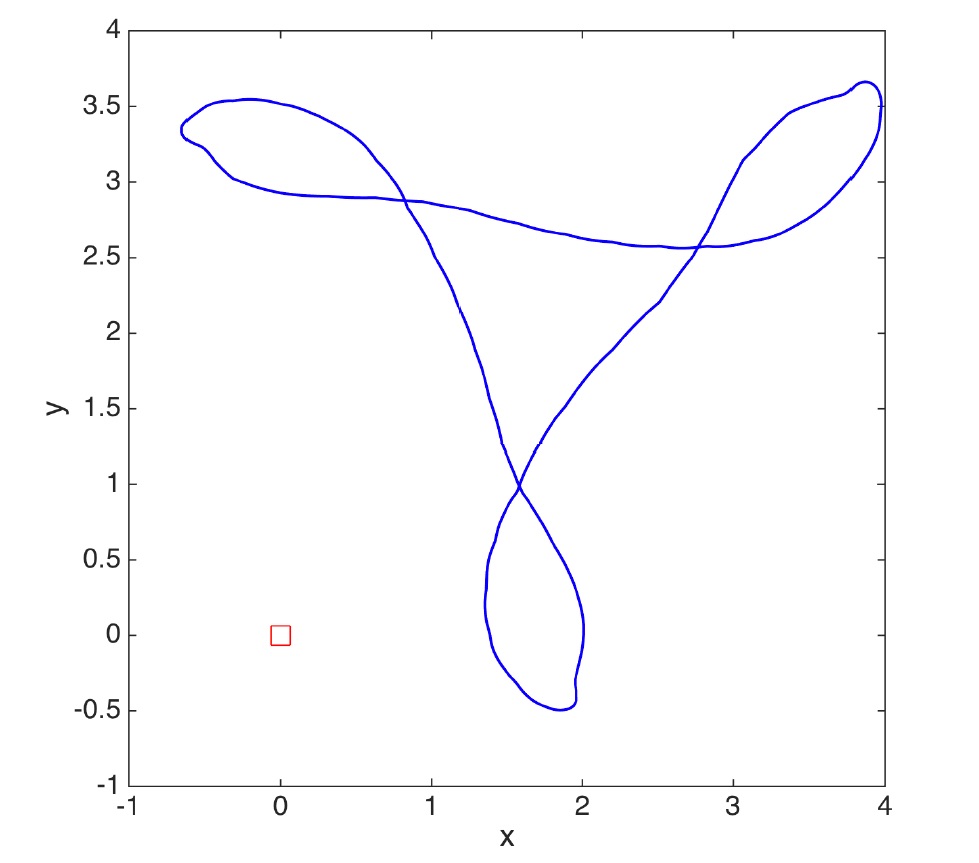}
\end{center}
\caption{Meander paths of spiral waves in numerical simulations of (\ref{FNPDE}) with inhomogeneity coefficients given by (\ref{coeffsinhsim1}).  The transients have been removed, and what is represented here is the final (steady) epicyclic meandering regime.  The 
kinetic parameters are $\beta=0.8$ and $\gamma=0.5$ for both, and $\tau=0.1018$ (left) and $\tau=0.10195$ (right).  The 2-frequency 3-petal motion thus appears to be phase-locked, as predicted by Theorem \ref{thmmainphaselock}. The red square at $(0,0)$ indicates a lattice point.}
\label{figpl1}
\end{figure}

In Figure \ref{figpl1}, we illustrate meander paths
for two simulations of (\ref{FNPDE}) with different parameters, and with transients removed. We observe a phase-locked 3-petal meandering path.  The inhomogeneity functions $g_1$ and $g_2$ in (\ref{gspecif})
have the following coefficients for both simulations:
\begin{equation}
\ep=0.01,\,\,A_1=-0.1997,\,\,A_2=0.2997,\,\,B_1=0.001,\,\,B_2=-0.001,\,\,C_1=-1,\,\,C_2=1.5.
\label{coeffsinhsim1}
\end{equation}

\subsection{Symmetry properties of (\ref{aveqsrat})}

The following is a direct consequence of the symmetry properties (\ref{symproplat}) for the functions $G_1$ and $G_2$ in (\ref{preav_nonint_rational}) and the definitions
(\ref{averages_rat}).
\begin{prop}
The averaged differential equations (\ref{aveqsrat}) are such that
\[
{\cal G}_1(-J\Psi,\varphi+\pi/2)=-J {\cal G}_1(\Psi,\varphi)\,\,\,\,\mbox{\rm and}\,\,\,\,
{\cal G}_2(-J\Psi,\varphi+\pi/2)={\cal G}_2(\Psi,\varphi)
\]
\label{propavsymrat}
\end{prop}
\proof This is a straightforward computation.
\hfill \qed

It follows from this proposition that if $(\Psi_0,\varphi_0)$ is an equilibrium for (\ref{aveqsrat}), then there are conjugate equilibria $(-J^n\Psi_0,\varphi_0+n\pi/2)$, $n=0,1,2,3$. Consequently,
the $2\pi\ell$-periodic solutions (\ref{solpredav}) of (\ref{preav_nonint_rational}) also come in conjugate families, and it follows from a straightforward computation that
these possess the symmetry properties
\begin{equation}
\begin{array}{c}
f^{1}_{(-J\Psi_0,\varphi_0+\pi/2)}(\theta,\ep)=-Jf^{1}_{(\Psi_0,\varphi_0)}(\theta,\ep),\\[0.15in]
f^{2}_{(-J\Psi_0,\varphi_0+\pi/2)}(\theta,\ep)=f^{2}_{(\Psi_0,\varphi_0)}(\theta,\ep)+\pi/2.
\end{array}
\label{solpredavsyms}
\end{equation}

Now, depending on the denominator $\ell$ for $\omega=k/\ell$, the functions ${\cal G}_{1}$ and ${\cal G}_2$ may possess additional symmetries.
\begin{prop}
If $\ell$ is an integer multiple of $4$, then 
\begin{equation}
{\cal G}_1(J\Psi,\varphi)=J {\cal G}_1(\Psi,\varphi)\,\,\,\,\mbox{\rm and}\,\,\,\,
{\cal G}_2(J\Psi,\varphi)={\cal G}_2(\Psi,\varphi).
\label{specsym1}
\end{equation}
If $\ell$ is an integer multiple of $2$, then
\begin{equation}
{\cal G}_1(-\Psi,\varphi)=- {\cal G}_1(\Psi,\varphi)\,\,\,\,\mbox{\rm and}\,\,\,\,
{\cal G}_2(-\Psi,\varphi)={\cal G}_2(\Psi,\varphi).
\label{specsym2}
\end{equation}
\end{prop}
\proof We will prove (\ref{specsym1}). The proof of (\ref{specsym2}) is similar, so we will omit it.  Let $n_1$ and $n_2$ be integers such that
\[
4\ell n_1 + 4k n_2 = \ell.
\]
Such integers $n_1$ and $n_2$ exist since $\mbox{\rm gcd}(4\ell,4k)=4$, and $\ell$ is a multiple of 4.  It follows that
\[
\begin{array}{rcl}
{\cal G}_1(J\Psi,\varphi)&=&{\displaystyle\frac{1}{2\pi\ell}\int_0^{2\pi\ell}\,R_{\varphi}R_{k\theta/\ell} G_1(J\Psi,\varphi+k\theta/\ell,\theta,0)\,d\theta}\\[0.15in]
&=&{\displaystyle\frac{1}{2\pi\ell}\int_0^{2\pi\ell}\,R_{\varphi}R_{k\theta/\ell} G_1(J\Psi,\varphi+(-2n_1+1/2)\pi+k\theta/\ell-\pi/2,\theta,0)\,d\theta}\\[0.15in]
&=&{\displaystyle\frac{1}{2\pi\ell}\int_0^{2\pi\ell}\,R_{\varphi}R_{k\theta/\ell} G_1(\Psi,\varphi+(-2n_1+1/2)\pi+k\theta/\ell,\theta,0)\,d\theta}\\[0.15in]
&=&{\displaystyle\frac{1}{2\pi\ell}\int_0^{2\pi\ell}\,R_{\varphi}R_{k\theta/\ell} G_1\left(\Psi,\varphi+\frac{k}{\ell}\left(\theta+\frac{\ell}{2k}(-4n_1+1)\pi\right),\theta,0\right)\,d\theta}\\[0.15in]
&=&{\displaystyle\frac{1}{2\pi\ell}\int_{\frac{\ell}{2k}(-4n_1+1)\pi}^{2\pi\ell+\frac{\ell}{2k}(-4n_1+1)\pi}\,
R_{\varphi}R_{k\theta/\ell}R_{2n_1\pi}R_{-\pi/2} G_1(\Psi,\varphi+k\theta/\ell,\theta+2\pi n_2,0)\,d\theta}\\[0.15in]
&=&{J\left(\displaystyle\frac{1}{2\pi\ell}\int_{2\pi n_2}^{2\pi\ell+2\pi n_2}\,R_{\varphi}R_{k\theta/\ell} G_1(\Psi,\varphi+k\theta/\ell,\theta,0)\,d\theta\right)}\\[0.15in]
&=&{J\left(\displaystyle\frac{1}{2\pi\ell}\int_0^{2\pi\ell}\,R_{\varphi}R_{k\theta/\ell} G_1(\Psi,\varphi+k\theta/\ell,\theta,0)\,d\theta\right)}=J{\cal G}_1(\Psi,\varphi).
\end{array}
\]
A similar computation using the second equation in (\ref{averages_rat}) reveals that ${\cal G}_2(J\Psi,\varphi)={\cal G}_2(\Psi,\varphi)$.
\hfill\qed

These additional symmetry properties have the following consequence:
\begin{cor}
If $\omega=k/\ell$, ($\mbox{\rm gcd}(k,\ell)=1$) is such that $\ell$ is even, then ${\cal G}_1$ in (\ref{aveqsrat}) is such that
\[
{\cal G}_1(0,\varphi)=0,\,\,\,\,\forall\varphi\in \Sone.
\]
Therefore, the roots of the $\pi/2$ periodic function ${\cal G}_2(0,\varphi)$ 
correspond to equilibria of (\ref{aveqsrat}). 
If $\varphi_0$ is such that ${\cal G}_2(0,\varphi_0)=0$, and the equilibrium point $(0,\varphi_0)$ for (\ref{aveqsrat}) is hyperbolic, then
(\ref{solpredavsyms}) reduce to
\begin{equation}
\begin{array}{c}
f^{1}_{(0,\varphi_0+\pi/2)}(\theta,\ep)=-Jf^{1}_{(0,\varphi_0)}(\theta,\ep),\\[0.15in]
f^{2}_{(0,\varphi_0+\pi/2)}(\theta,\ep)=f^{2}_{(0,\varphi_0)}(\theta,\ep)+\pi/2.
\end{array}
\label{solpredavsyms2}
\end{equation}
\end{cor}
\proof
This is a straightforward computation, using (\ref{specsym1}) and (\ref{specsym2}), and the fact that $J$ and $J^2=-I$ are rotation matrices which leave only the origin invariant.
\hfill\qed

\begin{prop}
Let $\omega=k/\ell$, ($\mbox{\rm gcd}(k,\ell)=1$) be such that $\ell$ is even, and let $(0,\varphi_0)$ be a hyperbolic equilibrium point for (\ref{aveqsrat}).  Then the
$2\pi\ell$-periodic solution (\ref{solpredav}) of (\ref{preav_nonint_rational})
has the spatio-temporal symmetry
\[
\begin{array}{c}
f^1_{(0,\varphi_0)}(\theta+\pi\ell,\ep)=-f^1_{(0,\varphi_0)}(\theta,\ep)\\[0.15in]
f^2_{(0,\varphi_0)}(\theta+\pi\ell,\ep)=f^2_{(0,\varphi_0)}(\theta,\ep).
\end{array}
\]
If $\ell$ is an integer multiple of 4, we have the additional spatio-temporal symmetry
\[
\begin{array}{c}
f^1_{(0,\varphi_0)}(\theta+\pi\ell/2,\ep)=(-J)^kf^1_{(0,\varphi_0)}(\theta,\ep)\\[0.15in]
f^2_{(0,\varphi_0)}(\theta+\pi\ell/2,\ep)=f^2_{(0,\varphi_0)}(\theta,\ep).
\end{array}
\]
\end{prop}
\proof
Again, we will only prove the case where $\ell$ is an integer multiple of 4, the other case being similar.  We note that the change of variables
$\theta=\tilde{\theta}+\pi\ell/2$, $\tilde{\varphi}=\varphi+k\pi/2$ leaves (\ref{preav_nonint_rational}) invariant.  From this, it immediately follows that
\[
f^1_{(0,\varphi_0)}(\theta,\ep)=f^1_{(0,\tilde{\varphi_0})}(\tilde{\theta},\ep),\,\,\,\,\,
f^2_{(0,\varphi_0)}(\theta,\ep)+k\pi/2=f^2_{(0,\tilde{\varphi_0})}(\tilde{\theta},\ep),
\]
or equivalently
\begin{equation}
f^1_{(0,\varphi_0)}(\theta+\pi\ell/2,\ep)=f^1_{(0,\varphi_0+k\pi/2)}(\theta,\ep),\,\,\,\,\,
f^2_{(0,\varphi_0)}(\theta+\pi\ell/2,\ep)+k\pi/2=f^2_{(0,\varphi_0+k\pi/2)}(\theta,\ep).
\label{interm1}
\end{equation}
But it follows from (\ref{solpredavsyms2}) that
\begin{equation}
\begin{array}{c}
f^{1}_{(0,\varphi_0+k\pi/2)}(\theta,\ep)=(-J)^kf^{1}_{(0,\varphi_0)}(\theta,\ep)\\[0.15in]
f^{2}_{(0,\varphi_0+k\pi/2)}(\theta,\ep)=f^{2}_{(0,\varphi_0)}(\theta,\ep)+k\pi/2.
\end{array}
\label{interm2}
\end{equation}
The conclusion follows from combining (\ref{interm1}) with (\ref{interm2}).
\hfill\qed

\noindent
{\bf Interpretation}:
For spiral waves, lattice symmetry-breaking can lead to phase-locked meandering spiral waves with even numbered ``petals'' which are anchored at a lattice point, and have spatio-temporal symmetry
characterized as follows: rotating the meander path about the lattice point by an angle of $\pi$ (or $\pi/2$ if the number of petals is a multiple of 4) is the same as advancing in time along the meander path by half (or a quarter) of the period.  It follows that the meander paths as a whole are invariant under rotations by $\pi$ (or $\pi/2$ if the number of
petals is a multiple of 4).

\begin{figure}[htpb]
\begin{center}
\includegraphics[width=3.2in]{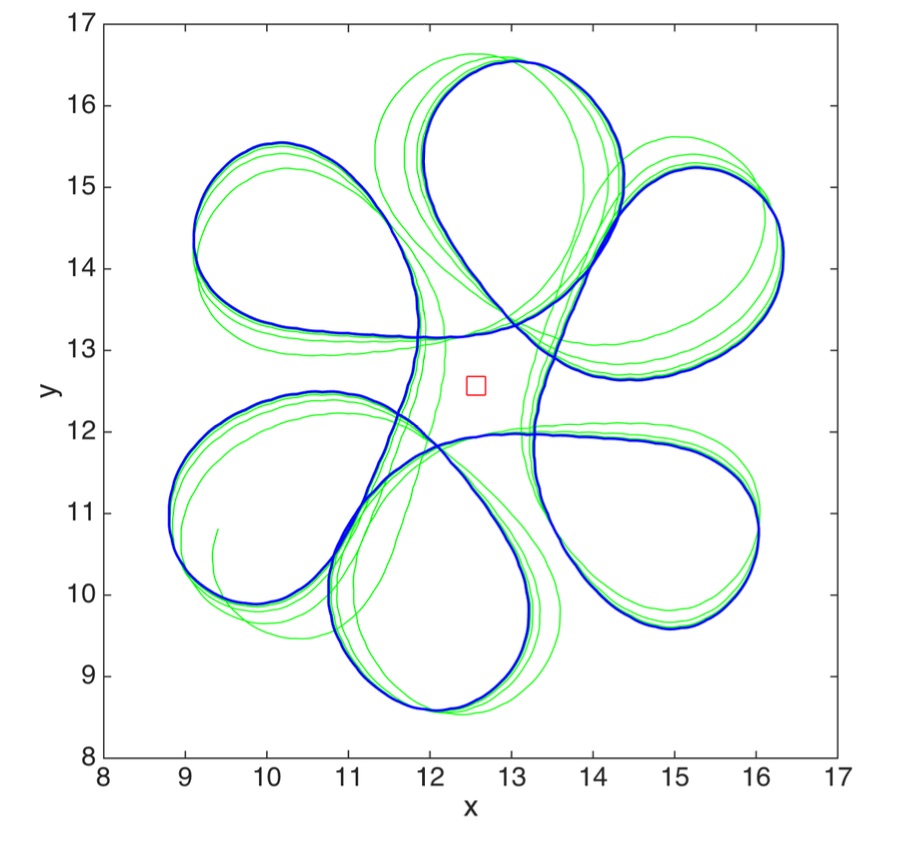}
\includegraphics[width=3.2in]{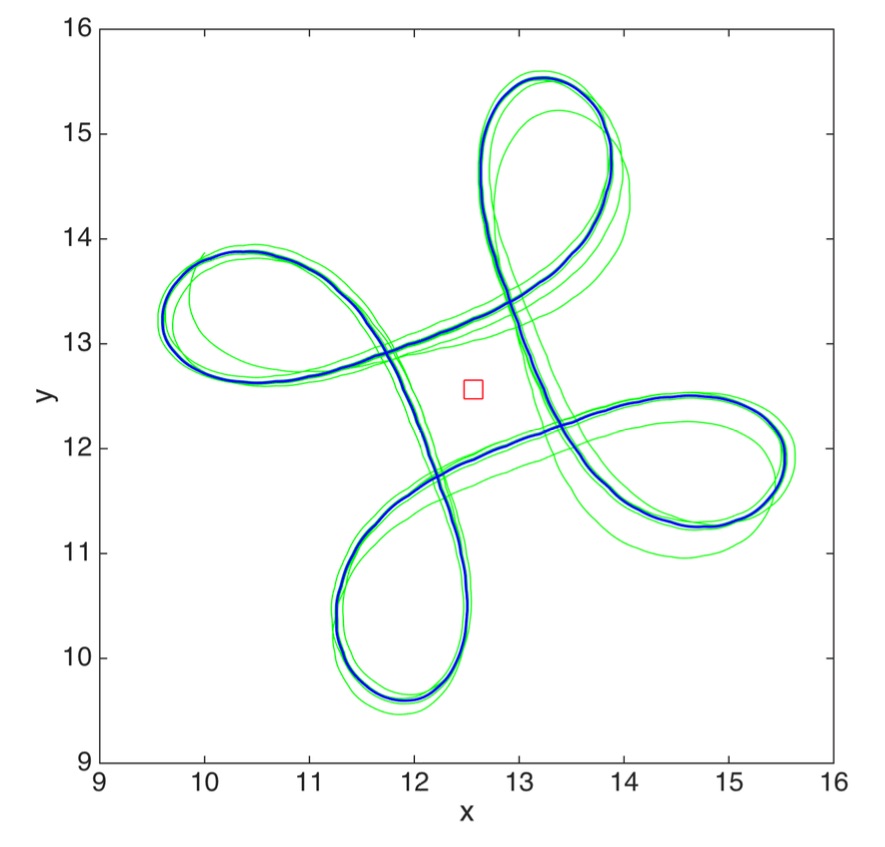}
\end{center}
\caption{Meander paths of spiral waves in numerical simulations of (\ref{FNPDE}) with inhomogeneity coefficients given by (\ref{coeffsinhsim2}) (left) and by (\ref{coeffsinhsim3}) (right).  The transients are in green and the final phase-locked meander path is in blue.  The 
kinetic parameters are $(\beta,\gamma,\tau)=(0.87,0.49,0.22)$ (left) and $(\beta,\gamma,\tau)=(0.8,0.65,0.17)$ (right).  The red square at $(4\pi,4\pi)$ indicates a lattice point.}
\label{figpl2}
\end{figure}

In Figure \ref{figpl2}, we illustrate meander paths
for two simulations of (\ref{FNPDE}) with different parameters and inhomogeneity functions. We observe 
in one case a phase-locked 6-petal meandering path, and in another case a phase-locked 4-petal meandering path.  The rotational symmetries of these meander paths are
consistent with the above remarks.  The inhomogeneity functions $g_1$ and $g_2$ in (\ref{gspecif}) are
\begin{equation}
\ep=0.01,\,\,A_1=1.4,\,\,A_2=0.92,\,\,B_1=5,\,\,B_2=1,\,\,C_1=-1,\,\,C_2=3.
\label{coeffsinhsim2}
\end{equation}
for the six-petal path, and
\begin{equation}
\ep=0.01,\,\,A_1=-1.75,\,\,A_2=-0.35,\,\,B_1=-6.25,\,\,B_2=1.25,\,\,C_1=1.25,\,\,C_2=-3.75.
\label{coeffsinhsim3}
\end{equation}
for the four-petal path.

\subsection{Phase-locking windows and bifurcations}

As we have seen above, linearly stable equilibrium points of (\ref{aveqsrat}) correspond to linearly stable phase-locked commensurate frequency meandering solutions of (\ref{cb-sys}).  We are now interested in understanding how these solutions behave under variation of $\omega$ and $\ep$.

In (\ref{preav_nonint}), we introduce an additional detuning parameter $\zeta$, and write $\omega=k/\ell+\ep\zeta$ (with $\mbox{\rm gcd}(k,\ell)=1$).  By following the same
procedure as we did at the beginning of this section, we obtain a parametrized system of averaged equations analogous to (\ref{aveqsrat}):
\begin{equation}
\begin{array}{rcl}
\dot{\Psi}&=&\ep {\mathscr G}_1(\Psi,\varphi,\zeta)\\
\dot{\varphi}&=&\ep {\mathscr G}_2(\Psi,\varphi,\zeta).
\end{array}
\label{aveqsrat_pert}
\end{equation}
where the functions ${\mathscr G}_{1,2}$ possess the same symmetry properties as the functions ${\cal G}_{1,2}$ (respectively), and
${\mathscr G}_{1,2}(\Psi,\varphi,0)={\cal G}_{1,2}(\Psi,\varphi)$.

Therefore,
if $(\omega,\ep)=(k/\ell,\ep_0)$ (with $\ep_0>0$, $\mbox{\rm gcd}(k,\ell)=1$) is such that (\ref{aveqsrat}) possesses such a linearly stable equilibrium point $(\Psi_0,\varphi_0)$,
then by using the implicit function theorem, there is a neighbourhood of the point $(0,\ep_0)$ in $(\zeta,\ep)$-space where (\ref{aveqsrat_pert}) possesses a linearly
stable equilibrium point $(\Psi(\zeta,\ep),\varphi(\zeta,\ep))$ close to $(\Psi_0,\varphi_0)$.  From this argument, it follows that the $\ell$-petal phase-locked meandering solution (\ref{workback}) of (\ref{cb-sys}) which corresponds to
$(\Psi_0,\varphi_0)$ persists in a neighbourhood of $(k/\ell,\ep_0)$ in the $(\omega,\ep)$ parameter space of (\ref{cb-sys})
We call such a neighbourhood a
{\em phase-locking window} for (\ref{cb-sys}).  The boundary of such a phase-locking window thus corresponds to bifurcation of the equilibrium point $(\Psi(\zeta,\ep),\varphi(\zeta,\ep))$ of (\ref{aveqsrat_pert}).  Since the phase space of (\ref{aveqsrat_pert}) is three-dimensional, then we could potentially observe higher-codimension
bifurcations (e.g. mode interactions) in (\ref{aveqsrat_pert}), and a thorough analysis will likely require the theoretical tools of \cite{LM} and \cite{WLM}.  We will not pursue this issue further in this paper.  Instead, we will briefly address saddle-node and Hopf
bifurcations.

\vspace*{0.15in}
\noindent
{\bf Saddle-node bifurcation}
\vspace*{0.15in}

\begin{figure}[htpb]
\begin{center}
\includegraphics[width=3.2in]{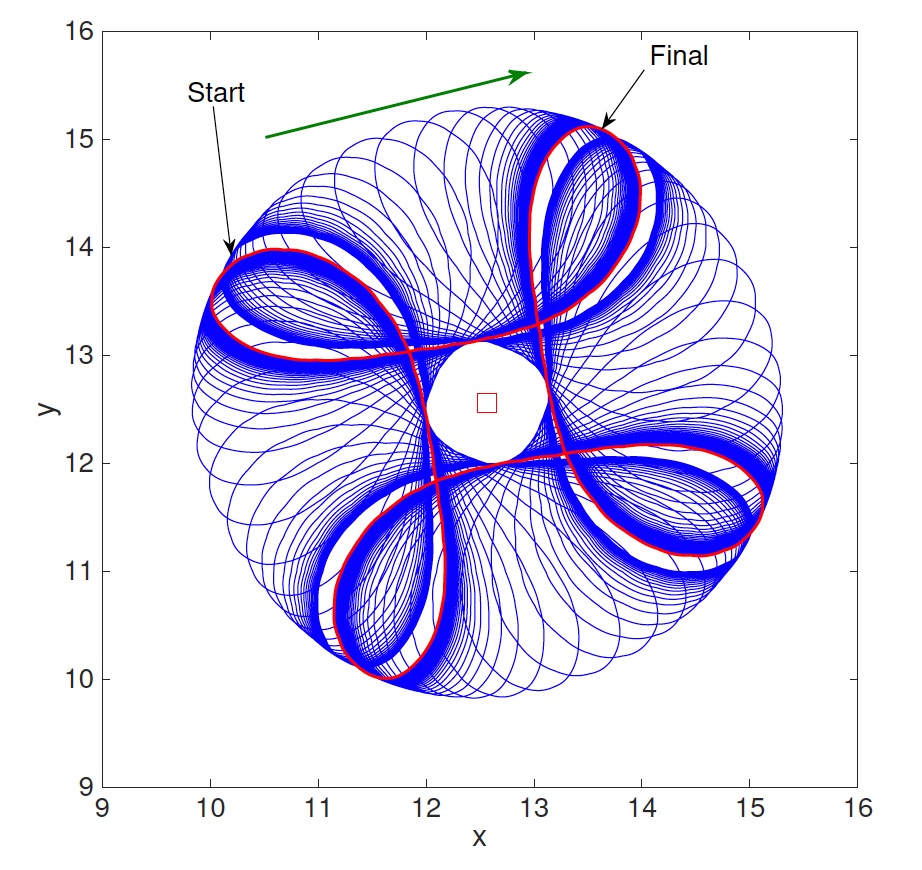}
\includegraphics[width=3.2in]{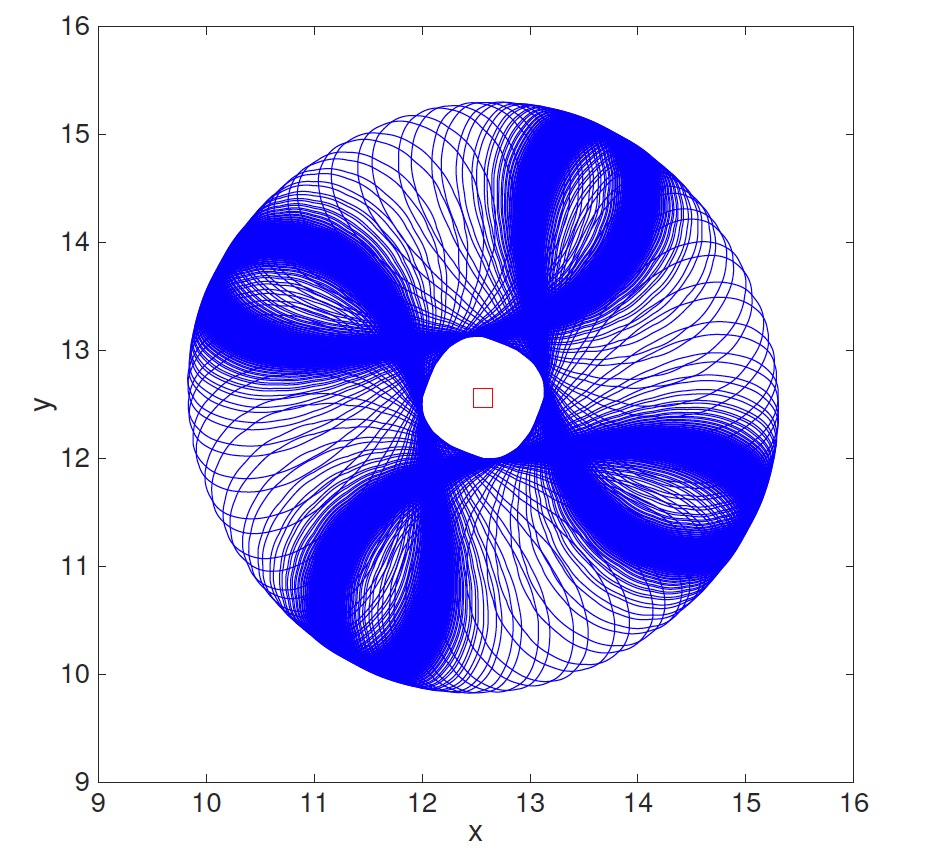}
\end{center}
\caption{Saddle-node bifurcation of 4-petal phase-locked meander paths in (\ref{FNPDE}) with inhomogeneity coefficients (\ref{coeffsinhsim5}).  For both simulations,
we have $\beta=0.8$ and $\gamma=0.65$.  In the left figure ($\tau=0.15818$), the meandering starts near an unstable 4-petal flower (presumably located in the white gap), flows in the direction of the green arrow, and eventually settles onto the stable 4-petal meandering state in red.  In the right figure ($\tau=0.15816$), the unstable and stable 4-petal flowers have coalesced and disappeared in a 
saddle-node bifurcation, and the meandering is now unlocked (although we can still
observe the slow transient passing through the remnants of the bifurcated states).}
\label{SN1}
\end{figure}

Phase-locked meander paths may lose stability via a saddle-node bifurcation of equilibria in (\ref{aveqsrat_pert}).  This phenomenon is very similar to the case where phase-locking is lost after exiting an Arnol'd tongue \cite{Arnold} in the theory of forced oscillators, so we will not pursue it in more details here.  In Figure \ref{SN1}, we see such a
saddle-node bifurcation in (\ref{FNPDE}) with inhomogeneity coefficients (\ref{gspecif}) given by (\ref{coeffsinhsim5}).

\vspace*{0.15in}
\noindent
{\bf Hopf bifurcation}
\vspace*{0.15in}

Here, we present what we believe to be a previously undocumented state for spiral waves.
Suppose that for $\ep=\ep_0>0$ and $\zeta=\zeta_0$, the point $(\Psi_0,\varphi_0)$ is an equilibrium point of (\ref{aveqsrat_pert}), and that the linearization of 
(\ref{aveqsrat_pert}) at $(\Psi_0,\varphi_0,\zeta_0)$ has a pair of complex conjugate eigenvalues on the imaginary axis and the other eigenvalue with negative real part.
Then generically, as $\zeta$ varies near $\zeta_0$, there will be a Hopf bifurcation, which will generate
a limit cycle with amplitude of the order $O(\sqrt{|\zeta-\zeta_0|})$.  

Heuristically, as a leading-order approximation, we write this limit cycle as
\[
\Psi(t)=\Psi_0+\sqrt{|\zeta-\zeta_0|}\,{\mathscr A}(t),\,\,\,\,\,\,
\varphi(t)=\varphi_0+\sqrt{|\zeta-\zeta_0|}\,{\mathscr B}(t)
\]
where ${\mathscr A}$ and ${\mathscr B}$ are ${\mathscr T}$-periodic, which leads to a leading-order approximation for the corresponding solution of (\ref{cb-sys}) 
\begin{equation}
\Psi_{\mbox{\tiny HB}}(t)=R_{\sqrt{|\zeta-\zeta_0|}\,{\mathscr B}(t)}\Psi(t),
\label{HBmeander}
\end{equation}
where $\Psi(t)$ is the $\ell$-petal phase-locked path given by (\ref{workback}).
In general, we expect that the period ${\mathscr T}$ will be incommensurate with $2\pi\ell$, and the resulting meander path (\ref{HBmeander}) can best be described as a ``fattened'' $\ell$-petal flower, with thickness $O(\sqrt{|\zeta-\zeta_0|})$: this fattening results from the small-amplitude time-periodic rotation $R_{\sqrt{|\zeta-\zeta_0|}\,{\mathscr B}(t)}$
of the overall $\ell$-petal meander path (\ref{workback}) that is losing stability at the bifurcation.  In the vocabulary of bifurcation theory, (\ref{HBmeander}) represents a Naimark-Sacker bifurcation from the $2\pi\ell$-periodic phase-locked solution (\ref{workback}).

\begin{figure}[htpb]
\begin{center}
\includegraphics[width=4.0in]{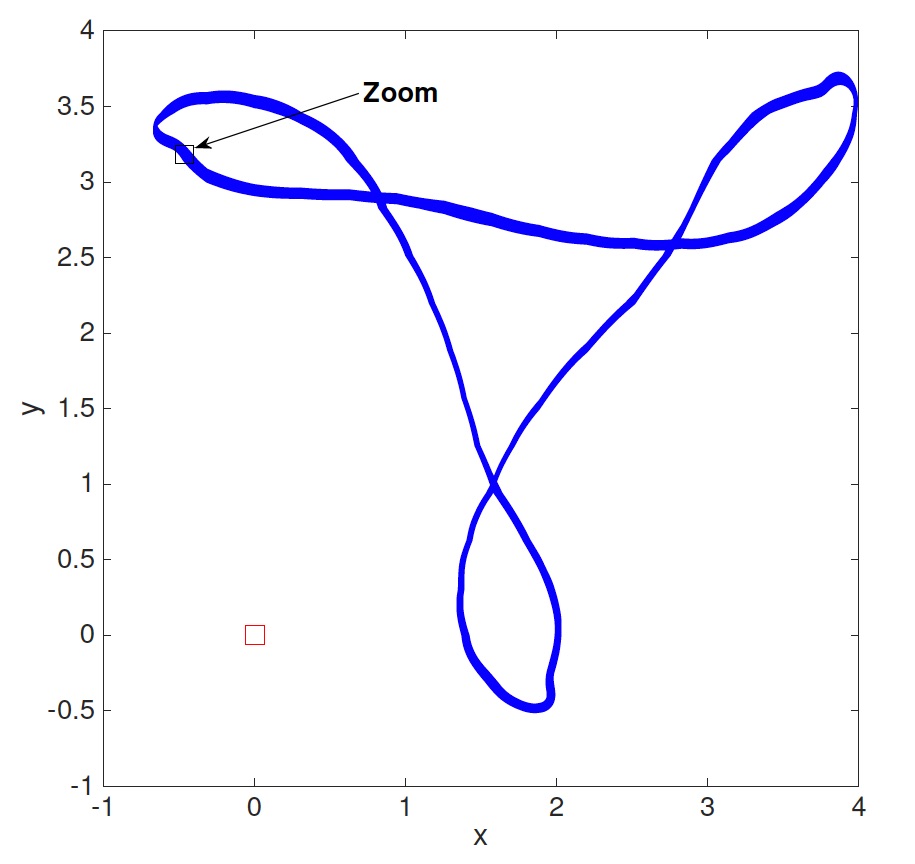}
\includegraphics[width=4.0in]{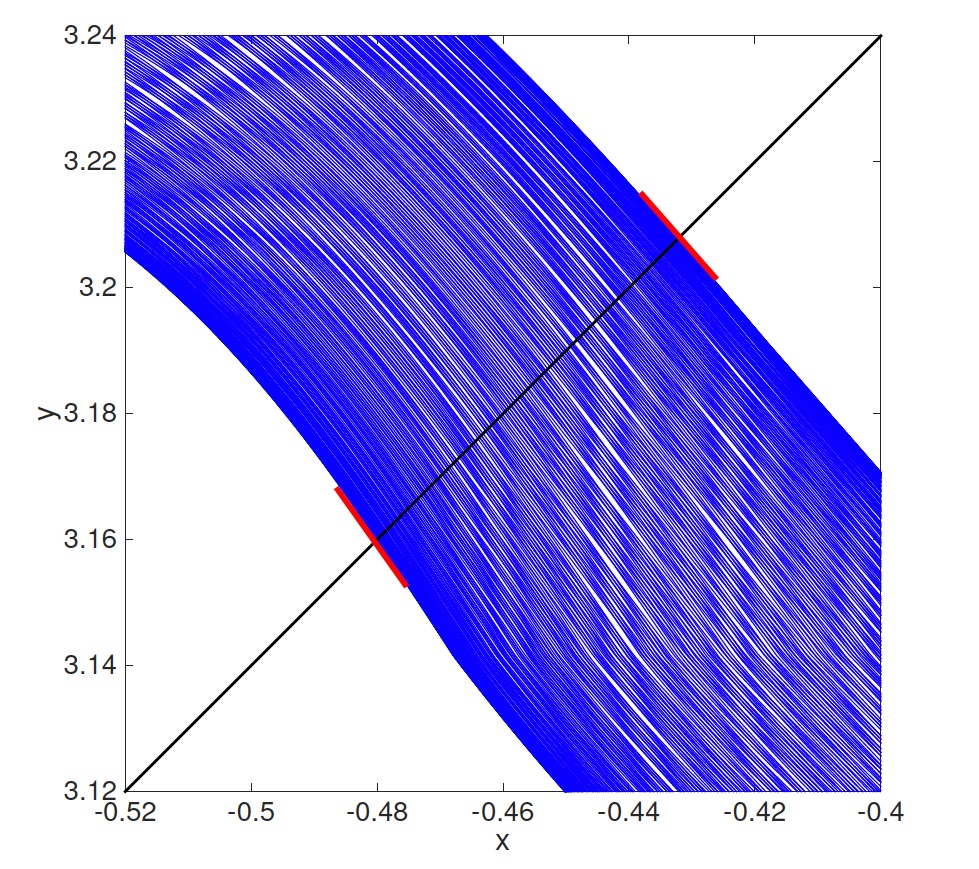}
\end{center}
\caption{Fattened 3-petal meander path observed in a simulation of (\ref{FNPDE}) with inhomogeneity (\ref{coeffsinhsim4}) and kinetic parameters
$\tau=0.1014$, $\beta=0.8$, $\gamma=0.5$ (transients removed).  The bottom figure is a close-up view of the inset box in the top figure.  As an admittedly arbitrary measure of the thickness, we choose
the length between the two little red lines of the black segment (line connecting the point $(-0.52,3.12)$ to the point $(-0.40,3.24)$).  In this case, it measures $0.06760$. In Figure
\ref{HB2} below, we illustrate how this meander path originates from a Hopf bifurcation from a phase-locked 3-petal meander path.}
\label{HB1}
\end{figure}

\begin{figure}[htpb]
\begin{center}
\includegraphics[width=3.2in]{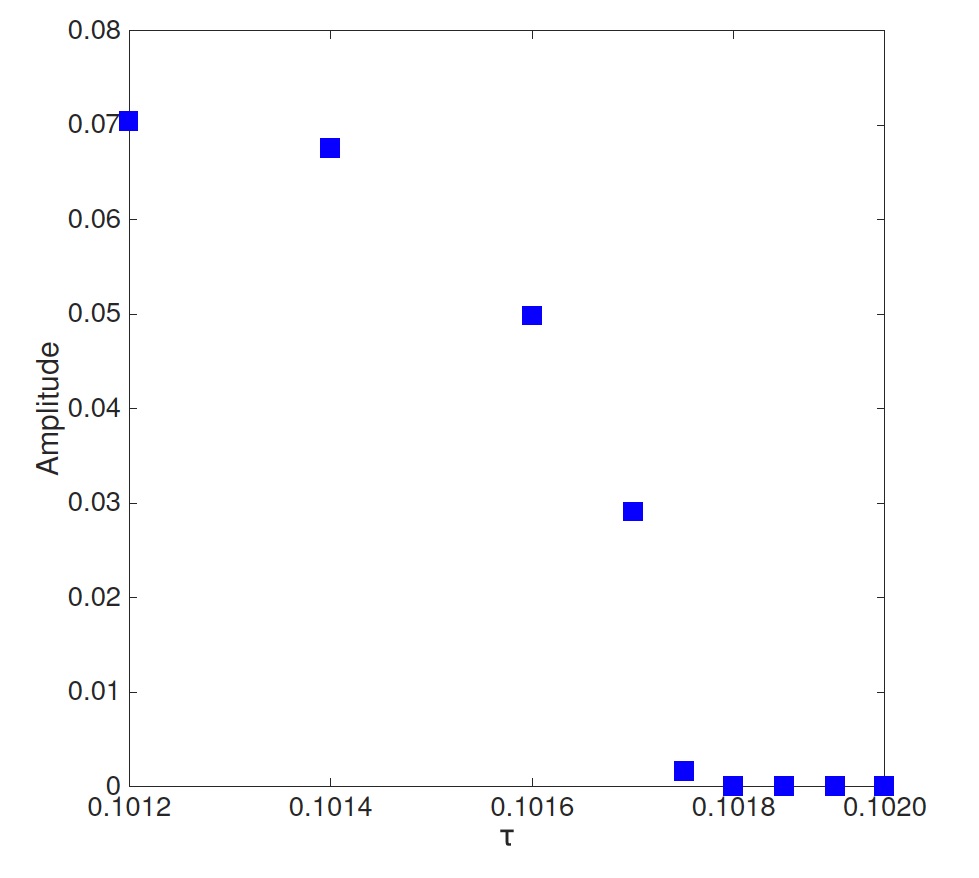}
\end{center}
\caption{Amplitude (measured as described in the caption of Figure \ref{HB1}) of fattened 3-petal meander path as a function of $\tau$ for 9 simulations of (\ref{FNPDE}).  The inhomogeneity data and values for $\beta$ and $\gamma$ are the same as in Figure \ref{HB1}, and $\tau$ is varied from $0.1012$ to $0.1020$.  We observe the characteristic shape of a Hopf bifurcation diagram.}
\label{HB2}
\end{figure}

\begin{figure}[htpb]
\begin{center}
\includegraphics[width=3.2in]{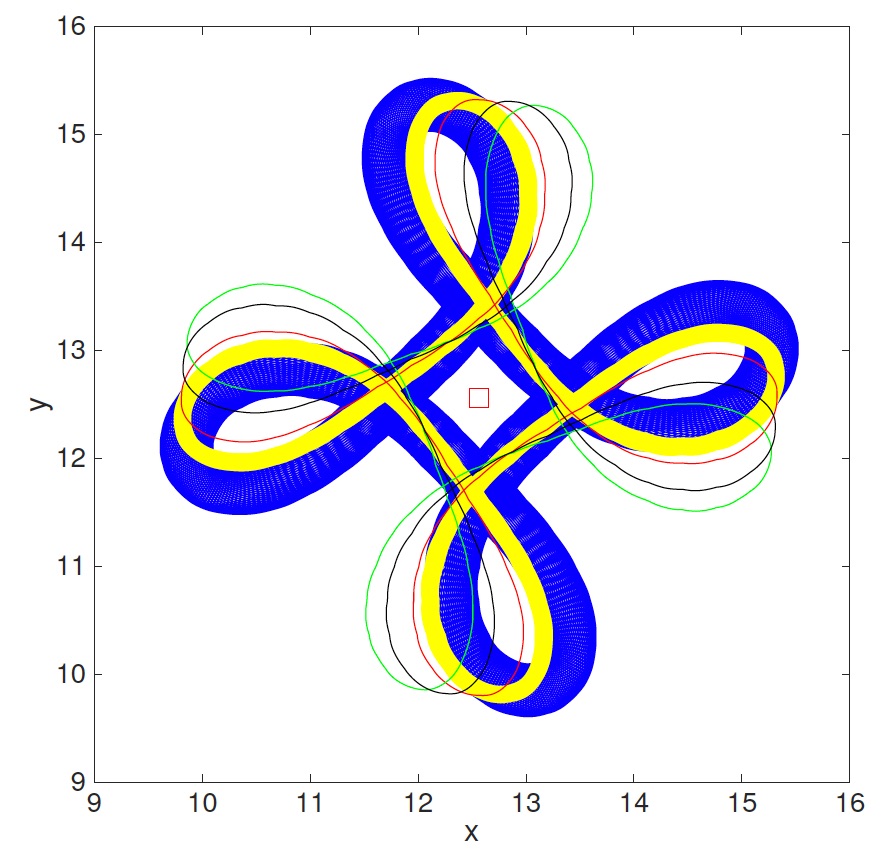}
\end{center}
\caption{Hopf bifurcation of a phase-locked symmetric 4-petal meander path in (\ref{FNPDE}) with inhomogeneity data (\ref{coeffsinhsim5}). For each of the five simulations,
we have $\beta=0.8$ and $\gamma=0.65$. The phase-locked paths (prior to bifurcation) are the green ($\tau=0.1586$), black ($\tau=0.1590$) and red ($\tau=0.1596$).
The bifurcated (fattened) paths are the yellow ($\tau=0.1600$) and blue ($\tau=0.1604$).}
\label{HB3}
\end{figure}

In Figure \ref{HB1}, we illustrate such a fattened $3$-petal flower in a simulation of (\ref{FNPDE}) with inhomogeneity coefficients (\ref{gspecif}) given by
\begin{equation}
\ep=0.01,\,\,A_1=-0.2,\,\,A_2=0.3,\,\,B_1=0,\,\,B_2=0,\,\,C_1=-1,\,\,C_2=1.5.
\label{coeffsinhsim4}
\end{equation}
In Figure \ref{HB2}, we illustrate how the image in Figure \ref{HB1} originates in a Hopf bifurcation from varying the parameter $\tau$ in (\ref{FNPDE}).
Figure \ref{HB3} illustrates a Hopf bifurcation from an anchored 4-petal flower for system (\ref{FNPDE}) with inhomogeneity data
\begin{equation}
\ep=0.01,\,\,A_1=-0.082,\,\,A_2=-0.014,\,\,B_1=-0.1,\,\,B_2=0.05,\,\,C_1=-0.25,\,\,C_2=-0.15
\label{coeffsinhsim5}
\end{equation}
in (\ref{gspecif}).


\Section{Stable linear meandering waves in (\ref{preav_int})}

When $\ep=0$ in (\ref{preav_int}), the solutions are such that 
\begin{equation}
\varphi(t)=\varphi(0),\,\,\theta(t)=t+\theta(0),\,\,\Psi(t)=R_{\varphi(0)}Vt+\Psi(0).
\label{3tflow}
\end{equation}
Working backwards through the changes of coordinates which transformed (\ref{cb-sys}) into (\ref{preav_int}), this corresponds to linearly meandering travelling waves, such as illustrated with the red meander path in Figure \ref{fig3}. In the context of full $\SEtwo$ symmetry (i.e. $\ep=0$), this
phenomenon is not structurally stable, since arbitrarily small changes to the value of $\omega$ in (\ref{cb-sys}) leads to a meander pattern which has (large radius) circular shape (green and blue meander paths in Figure \ref{fig3}).

In the sequel, it will be helpful to interpret (\ref{3tflow}) in the following manner.  Consider the hypersurface $\varphi=\mbox{\rm constant}=\varphi(0)$ in $\mathbb{T}^4$. This hypersurface is diffeomorphic to a three-torus (parametrized by $\psi_{1,2}$ and $\theta$). Equation (\ref{3tflow}) describes a linear flow on this three-torus, with 
frequencies given by the two components of $R_{\varphi(0)}V$ and $1$.  If these frequencies are rationally independent, then the solution curves for this linear flow are dense on the three-torus.  In a sense, we are interested in the persistence of such densely filled invariant three-tori when $\ep$ becomes non-zero.

Our main result in this section is the following:

\begin{theo}
Consider the system (\ref{preav_int}), and define the function
\[
{\cal Z}(\varphi)=\langle H_2(\cdot,\varphi,\cdot,0) \rangle \equiv \frac{1}{(2\pi)^3}\int_{\mathbb{T}^3}\,H_2(\Psi,\varphi,\theta,0)d\Psi d\theta
\]
Let $\varphi_0$ be such that ${\cal Z}(\varphi_0)=0$, and $\alpha={\cal Z}'(\varphi_0)\neq 0$.  Define the real numbers $\omega_1$ and $\omega_2$ as the components of the two-dimensional vector $R_{\varphi_0}V$ and suppose that $\boldsymbol{\Omega}=(\omega_1,\omega_2,1)$ satisfies a Diophantine condition such as in Definition \ref{DiophantineDef}.
Then for all $\ep>0$ sufficiently small, (\ref{preav_int}) has an invariant three-torus
\[
\varphi={\mathscr L}(\Psi,\theta,\ep),\,\,\,\,{\mathscr L}(\Psi,\theta,\ep)\rightarrow \varphi_0,\,\,\,\mbox{\rm as}\,\,\,\ep\rightarrow 0.
\]
This invariant three-torus is locally asymptotically stable (resp. unstable) if $\alpha < 0$ (resp. $\alpha > 0$).
\label{transverse}
\end{theo}
\proof
The proof follows closely that of Theorem 5.1 of \cite{CL}, with minor changes accounting for the additional angular variable $\theta$.  Therefore, we only give a sketch of the
main points.  First, we make the change of variables
\[
\varphi\rightarrow\varphi_0+\sqrt{\ep}\,\varphi
\]
and perform a Taylor expansion of (\ref{preav_int}) to get
\begin{equation}
\begin{array}{rcl}
\dot{\Psi}&=&R_{\varphi_0}V+\sqrt{\ep}R'_{\varphi_0}V\varphi+\ep (R_{\varphi_0}H_1(\Psi,\varphi_0,\theta,0)+\frac{1}{2}R''_{\varphi_0}V\varphi^2)+\ep^{\frac{3}{2}}Q_1(\Psi,\varphi,\theta,\ep)\\[0.1in]
\dot{\varphi}&=&\sqrt{\ep}H_2(\Psi,\varphi_0,\theta,0)+\ep H_{2,\varphi}(\Psi,\varphi_0,\theta,0)\varphi)+
\ep^{\frac{3}{2}}(\frac{1}{2}H_{2,\varphi\varphi}(\Psi,\varphi_0,\theta,0)\varphi^2+H_{2,\ep}(\Psi,\varphi_0,\theta,0))\\&&+\ep^2 S_1(\Psi,\varphi,\theta,\ep)\\[0.1in]
\dot{\theta}&=&1,
\end{array}
\label{preav_int_tayl}
\end{equation}
for smooth functions $Q_1$ and $S_1$.  Since the function $H_2(\Psi,\varphi_0,\theta,0)$ has zero mean value, i.e.  $\langle H(\cdot,\varphi_0,\cdot,0)\rangle = Z(\varphi_0)=0$, then it follows from Proposition 
\ref{DiophantineProp} that there exists a smooth function ${\cal Y}_1(\Psi,\theta)$ such that
\[
D_{\Psi}{\cal Y}_1(\Psi,\theta)R_{\varphi_0}V+D_{\theta}{\cal Y}_{1}(\Psi,\theta)=H_2(\Psi,\varphi_0,\theta,0).
\]
We now perform, in succession, the following changes of variables
\begin{itemize}
\item $\varphi\rightarrow \varphi+\sqrt{\ep}{\cal Y}_1(\Psi,\theta)$
\item $\varphi\rightarrow \sqrt{\ep}\varphi$
\item $\varphi\rightarrow \varphi+\ep ({\cal Y}_2(\Psi,\theta)+\varphi {\cal Y}_3(\Psi,\theta))$
\item $\Psi\rightarrow\Psi+\ep {\cal Y}_4(\Psi,\theta)$,
\end{itemize}
where ${\cal Y}_{2,3,4}$ are suitably chosen (via Proposition \ref{DiophantineProp}) to annihilate certain terms in the ODE. Equation (\ref{preav_int_tayl}) then becomes
\begin{equation}
\begin{array}{rcl}
\dot{\Psi}&=&R_{\varphi_0}V+\ep (R'_{\varphi_0}V\varphi+\kappa)+O(\ep^{\frac{3}{2}})\\
\dot{\varphi}&=&\ep (\alpha\varphi+c)+O(\ep^{\frac{3}{2}})\\
\dot{\theta}&=&1,
\end{array}
\label{preav_int_final}
\end{equation}
where $\alpha$ is as in the statement of the Theorem, and
\[
\kappa=\langle R_{\varphi_0}H_1(\cdot,\varphi_0,\cdot,0) \rangle,
\]
\[
c=\langle H_{2,\ep}(\cdot,\varphi_0,\cdot,0)+H_{2,\varphi}(\cdot,\varphi_0,\cdot,0){\cal Y}_1(\cdot,\cdot)
-D_{\Psi}{\cal Y}_1(\cdot,\cdot)(R_{\varphi_0}H_1(\cdot,\varphi_0,\cdot,0)+R'_{\varphi_0}V{\cal Y}_1(\cdot,\cdot) \rangle.
\]
A simple translation of the variable $\varphi$ in (\ref{preav_int_final}) renders the equation in the proper form to apply Theorem 2.3, \S VII.2 of \cite{Hale}  and get the conclusion.
\hfill\qed

\noindent
{\bf Interpretation}:  In fully $\SEtwo$ symmetric dynamical systems (such as (\ref{FNPDE}) with $\ep=0$), modulated travelling waves are not structurally stable.  This is illustrated in
figure 1 of \cite{BK}, where modulated travelling waves occur only on a line (dashed curve) in the $\beta-\tau$ parameter space of (\ref{FNPDE}) with $\ep=0$ and $\gamma=0.5$.
Ashwin, Melbourne and Nicol \cite{AMN} explain this as a codimension 1 drift bifurcation on group orbits of solutions to $\SEtwo$-equivariant dynamical systems.

\begin{figure}[htpb]
\begin{center}
\includegraphics[width=5.0in]{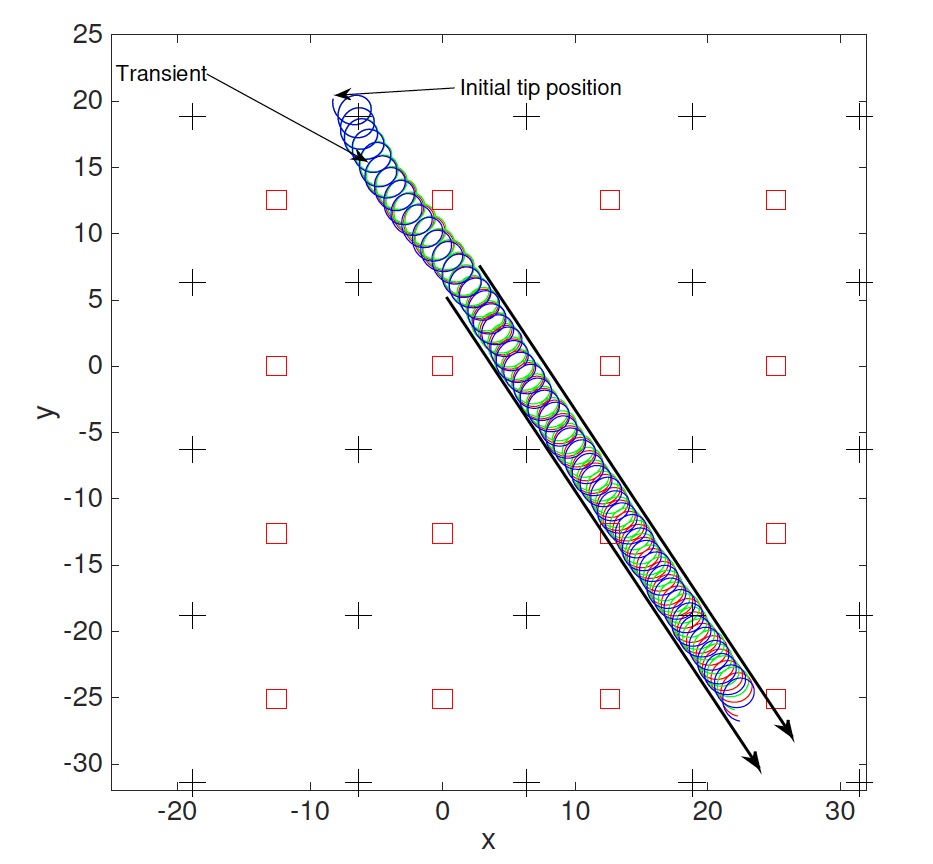}
\end{center}
\caption{Superposition of meander paths for three simulations of (\ref{FNPDE}) with 
inhomogeneity data given by (\ref{coefflin1})
 and same initial condition in all three cases. 
 These paths are consistent with a structurally stable modulated travelling wave, as predicted by Theorem \ref{transverse}. Parameter values are $\tau=0.26$, $\gamma=0.5$ and $\beta=0.793$ (green), $\beta= 0.79275$ (blue), $\beta=0.792875$ (red).
The parallel black arrows on both sides of the meander paths have been added to the figure to guide the eye, and the red squares and black crosses represent the lattice and dual lattice points, respectively.}
\label{linperts}
\end{figure}

Theorem \ref{transverse} implies that under certain conditions, it is possible for a lattice symmetry-breaking perturbation in (\ref{FNPDE}) to render modulated travelling waves structurally stable.
Figure \ref{linperts} appears to illustrate this phenomenon. In this figure, we show the meander paths for three numerical simulations of (\ref{FNPDE}) using the same kinetic parameters as
those of Figure \ref{fig3} (homogeneous case), but with inhomogeneity data (\ref{gspecif}) given by
\begin{equation}
\ep=0.01,\,\,A_1=-0.6,\,\,A_2=-0.4,\,\,B_1=-0.00005,\,\,B_2=0.00007,\,\,C_1=-3,\,\,C_2=-2
\label{coefflin1}
\end{equation}
and same initial condition for all three simulations.  This figure should be compared and contrasted to Figure \ref{fig3}.

\Section{Discussion} 

Discrete spatial structures are arguably the fundamental building blocks of nature.  In electrophysiological tissue, these structures are arrays of excitable cells coupled together through
gap junctions. In many instances, treating these structures as a homogeneous continuum, and deriving partial differential equations based on this assumption (using homogenization techniques) is a reasonable low-order approximation, and leads to many successful predictions about the system being modelled.  However, if the solutions we are interested in studying display features that are of a size comparable in order
to the cellular structure, then the continuum assumption may not be of sufficient complexity to capture the dynamical properties of these solutions.  Furthermore, numerical simulations of
spiral waves in homogeneous
RDPDEs using coarse spatial discretization may introduce spurious dynamical properties, which are inconsistent with a homogeneous model.
This paper is an attempt to characterize certain generic, qualitative phenomenological effects that a lattice structure may induce on meandering spiral waves in two-dimensional excitable media.
As was the case in previous studies \cite{BLM,CL,LeBlanc,LW}, we have shown that a weak lattice perturbation has stabilizing (in physical space) effects on the meandering motion, alters the
spatio-temporal symmetry properties of the meander path, and can lead to phase-locked solutions.  We also presented what we believe to be a previously undocumented case for spiral wave
meandering: the so-called ``fattened'' $\ell$-petal flower which originates in a Hopf bifurcation from a phase-locked $\ell$-petal epicyclic meandering wave.

As we mentioned in the Introduction, one could make a compelling case that perhaps different geometries of lattices would be more appropriate in certain arrangements of excitable cells
(e.g. hexagonal). However, 
the purpose of this paper was not to describe any one physical situation in particular, but to gain an understanding of how spatially organized discrete structures may affect spiral wave dynamics.
Certainly, one could repeat the analysis of this paper, replacing the square lattice by a hexagonal lattice.  We expect that the number theoretical aspects of which meander paths get anchored
at lattice points, as well as spatio-temporal symmetries would change, but we don't expect any fundamentally new qualitative results.  One could also adopt a purely discrete approach, for example
studying spirals in lattice dynamical systems (LDS) \cite{PE}.  But then one would lose the aspect that much of the observed spiral wave dynamics are driven by and originate in
Euclidean symmetry.  This is one
of the advantages of our forced symmetry-breaking approach: we are in a sense studying systems that interpolate the fully Euclidean continuum RDPDE, and completely discrete LDS.  
We thus expect that the results presented herein will be important in the as of yet largely unexplored study of meandering spiral waves in two-dimensional LDS.

We also note that the reduction of the infinite-dimensional problem (\ref{ade}) into the finite-dimensional center-bundle equations (\ref{cb-sys}) requires a spectral gap condition in 
Hypothesis \ref{hyp2}.  While this hypothesis holds for a large variety of spirals (such as decaying amplitude spirals), there is also a large family of spirals for which they don't, including
Archimedean spirals.  See \cite{Scheel} for a thorough discussion on these spectral issues.  Even if the spiral does not decay (e.g. Archimedean spirals), the
center bundle ODEs apparently still describe remarkably well many of the experimentally observed
dynamics of the wave.  Thus, the study of Euclidean-equivariant ODEs on
finite-dimensional non-compact manifolds, and their perturbations, is a central idea in the study of the dynamics
and bifurcations of spiral waves, and it is the one we have adopted here.

As a final remark, we note that one of the strengths of our approach in this paper is also perhaps one of its weaknesses, depending on one's purpose.  By adopting a model-independent
approach as we have done here (i.e. computations and analysis driven mostly by symmetry considerations), we capture some robust generic features of the effects of lattice symmetry-breaking for a
large class of mathematical models, without much consideration for the finer details of the models. The most important 
requirement is that we have a mathematical model where spiral waves are known to exist, and the model is weakly heterogeneous in space, with heterogeneities on a lattice.
So in this sense, our results are potentially broad in applicability.  However, one can imagine situations where it may be important, when studying a specific model (e.g. RDPDE such as (\ref{FNPDE})),
to know the explicit link between
the model parameters (e.g. $\tau$, $\beta$, $\gamma$, $\ep$, $g_1$ and $g_2$) and the quantities which arise in the center-bundle equations (\ref{cb-sys}).  As mentioned earlier in this paper
in Remark \ref{remark24}(a), although
such a link theoretically exists, it is in practice usually unfeasible to explicitly compute this link, unless one has explicit knowledge of the mathematical representation (e.g. algebraic formula) for the spiral wave.  

\vspace*{0.25in}
\noindent
{\Large\bf Acknowledgments}

\vspace*{0.2in}
This research is partly supported by the
Natural Sciences and Engineering Research Council of Canada in the
form of a Discovery Grant (VGL). One of the authors (PK) thanks the Department of Mathematics and Statistics at the University of Ottawa for hosting his postdoctoral fellowship during which this paper was written.

{\Large\bf\flushleft Appendix}

\appendix

\Section{Proofs of Propositions \ref{propnonint} and \ref{propint}}

The proofs for both these Propositions follow one another closely except for some subtle technical differences.  We will first give the proof of Proposition \ref{propnonint}, and then indicate
the modifications which need to be made in order to prove Proposition \ref{propint}.

The function $h_2(\theta)$ in (\ref{cb-sys}) satisfies (\ref{zm}). We recall the definition of the rotation matrix $R_{\phi}$ in (\ref{rotmatdef}), and that $J=R_{-\pi/2}$. 
We note that ${\displaystyle \frac{d}{d\phi} R_{\phi}=-JR_{\phi}}$.

If we define 
\begin{equation}
K(\theta)=R_{\int_0^{\theta}\,h_2(s)ds}h_1(\theta),
\label{KDEF}
\end{equation}
then it is easy to show that $K(\theta)$ is $2\pi$-periodic in $\theta$.  We then have
\begin{lemma}
If $\omega$ is not an integer, then there exists a $2\pi$-periodic function $M(\theta)$ such that
\begin{equation}
\frac{d}{d\theta}M(\theta)-\omega JM(\theta) = K(\theta).
\label{odeM}
\end{equation}
\end{lemma}
\proof
Write $K(\theta)$ and $M(\theta)$ as Fourier series ${\displaystyle K(\theta)=\sum_{k\in\mathbb{Z}}\,R_{k\theta}\hat{B}_k}$, 
${\displaystyle M(\theta)=\sum_{k\in\mathbb{Z}}\,R_{k\theta}\hat{A}_k}$, where the coefficients $\hat{A}_k$ and $\hat{B}_k$ are in $\mathbb{R}^2$.  Then (\ref{odeM}) is satisfied provided
$\hat{A}_k$ and $\hat{B}_k$ satisfy the algebraic equation
\[
-(k+\omega)J\hat{A}_k=\hat{B}_k\Longrightarrow \hat{A}_k=\frac{1}{k+\omega}J\hat{B}_k.
\]
Uniform convergence of the series for $M$ 
follows from the uniform convergence of the series for $K$.
\hfill\qed
\begin{lemma}
Let $M(\theta)$ be as in Lemma \ref{odeM}.  Define ${\cal S}(\theta)=R_{-\int_0^{\theta}\,h_2(s)ds}M(\theta)$
and ${\cal Q}(\varphi,\theta)=R_{\varphi}{\cal S}(\theta)$. Then ${\cal Q}$ is $2\pi$-periodic in both $\varphi$ and $\theta$, and satisfies the linear partial differential equation
\[
(\omega+h_2(\theta))\frac{\partial {\cal Q}}{\partial\varphi}+\frac{\partial {\cal Q}}{\partial\theta}=R_{\varphi}h_1(\theta).
\]
\label{lemQ}
\end{lemma}
\proof
The periodicity properties are obvious.
Using the fact that $M$ satisfies (\ref{odeM}) and $(R_{\phi})'=-JR_{\phi}$, we compute
\[
\frac{\partial {\cal S}}{\partial\theta} = h_1(\theta)+(\omega+h_2(\theta))J{\cal S}(\theta),
\]
from which it follows that
\[
\frac{\partial {\cal Q}}{\partial \varphi}=-JR_{\varphi}{\cal S}(\theta)\,\,\,\,\,\mbox{\rm and}\,\,\,\,\,
\frac{\partial {\cal Q}}{\partial \theta}=R_{\varphi} [h_1(\theta)+(\omega+h_2(\theta))J{\cal S}(\theta)].
\]
The conclusion now follows from a simple computation.
\hfill\qed

Now, setting 
\begin{equation}
\hat{\Psi}=\Psi-{\cal Q}(\varphi,\theta) (\mbox{\rm modulo}\,\, \mathbb{T}^2), 
\label{fundcoc}
\end{equation}
and using the above lemmas, the system (\ref{cb-sys}) becomes (upon dropping the hats)
\[
\begin{array}{rcl}
\dot{\Psi}&=&\ep R_{\varphi} {\mathscr M}_1(\Psi,\varphi,\theta,\ep)\\
\dot{\varphi}&=&\omega+h_2(\theta)+\ep {\mathscr M}_2(\Psi,\varphi,\theta,\ep)\\
\dot{\theta}&=&1,
\end{array}
\]
where
\[
{\mathscr M}_1(\Psi,\varphi,\theta,\ep)=F_1(\Psi+{\cal Q}(\varphi,\theta),\varphi,\theta,\ep)+J{\cal S}(\theta)F_2(\Psi+{\cal Q}(\varphi,\theta),\varphi,\theta,\ep)
\]
and
\[
{\mathscr M}_2(\Psi,\varphi,\theta,\ep)=F_2(\Psi+{\cal Q}(\varphi,\theta),\varphi,\theta,\ep).
\]
The functions ${\mathscr M}_{1,2}$ satisfy the symmetry properties (\ref{symproplat}) since the functions $F_{1,2}$ satisfy these properties, and
\[
{\cal Q}(\varphi+\pi/2,\theta)=R_{\varphi+\pi/2}{\cal S}(\theta)=R_{\pi/2}R_{\varphi}{\cal S}(\theta)=-J{\cal Q}(\varphi,\theta).
\]
Finally, setting $\hat{\varphi}=\varphi-\int_0^{\theta}\,h_2(s)ds$ (modulo $\mathbb{T}^1$) and dropping the hat, we obtain equations (\ref{propnonint}), where
\[
\begin{array}{c}
G_1(\Psi,\varphi,\theta,\ep)=R_{\int_0^{\theta}\,h_2(s)ds} {\mathscr M}_1(\Psi,\varphi+\int_0^{\theta}\,h_2(s)ds,\theta,\ep),\\[0.2in]
G_{2}(\Psi,\varphi,\theta,\ep)={\mathscr M}_{2}(\Psi,\varphi+\int_0^{\theta}\,h_2(s)ds,\theta,\ep)
\end{array}
\]
also satisfy the symmetry properties (\ref{symproplat}).
This ends the proof of Proposition \ref{propnonint}
\hfill\qed

To prove Proposition \ref{propint}, we define $K(\theta)$ as in (\ref{KDEF}).  Let $j\in\mathbb{Z}$ be such that $\omega=-j$.  If the Fourier series for $K(\theta)$ is
${\displaystyle K(\theta)=\sum_{k\in\mathbb{Z}}\,R_{k\theta}\hat{B}_k}$, then define
\[
\tilde{K}(\theta)=\sum_{k\in\mathbb{Z},k\neq j}\,R_{k\theta}\hat{B}_k
\]
so that
\[
K(\theta)=\tilde{K}(\theta)+R_{j\theta}\hat{B}_{j}.
\]
We then have
\begin{lemma}
If $\omega=-j\in\mathbb{Z}$, then there exists a $2\pi$-periodic function $\tilde{M}(\theta)$ such that
\[
\frac{d}{d\theta}\tilde{M}(\theta)-\omega J\tilde{M}(\theta) = \tilde{K}(\theta),
\]
\label{odeM2}
\end{lemma}
whose proof follows exactly like that of Lemma \ref{odeM}, considering that the Fourier series for $\tilde{K}$ is such that the term for $k=j$ vanishes, so we can set $A_{j}=0$ in the Fourier series for $\tilde{M}$. 

\medskip
Lemma \ref{lemQ} is replaced by
\begin{lemma}
Let $\tilde{M}(\theta)$ be as in Lemma \ref{odeM2}.  Define $\tilde{{\cal S}}(\theta)=R_{-\int_0^{\theta}\,h_2(s)ds}\tilde{M}(\theta)$
and $\tilde{{\cal Q}}(\varphi,\theta)=R_{\varphi}\tilde{{\cal S}}(\theta)$. Then $\tilde{{\cal Q}}$ is $2\pi$-periodic in both $\varphi$ and $\theta$, and satisfies the linear partial differential equation
\[
(\omega+h_2(\theta))\frac{\partial \tilde{{\cal Q}}}{\partial\varphi}+\frac{\partial \tilde{{\cal Q}}}{\partial\theta}=R_{\varphi}[h_1(\theta)-R_{-\int_0^{\theta}\,h_2(s)ds}R_{j\theta}B_j]
\]
\label{lemQ2}
\end{lemma}

Performing the change of variables $\hat{\Psi}=\Psi-\tilde{{\cal Q}}(\varphi,\theta)$ on (\ref{cb-sys}) yields (upon dropping the hats)
\[
\begin{array}{rcl}
\dot{\Psi}&=&R_{\varphi}R_{j\theta-\int_0^{\theta}\,h_2(s)ds}B_j+\ep R_{\varphi} {\mathscr N}_1(\Psi,\varphi,\theta,\ep)\\
\dot{\varphi}&=&\omega+h_2(\theta)+\ep {\mathscr N}_2(\Psi,\varphi,\theta,\ep)\\
\dot{\theta}&=&1,
\end{array}
\]
where
\[
{\mathscr N}_1(\Psi,\varphi,\theta,\ep)=F_1(\Psi+\tilde{{\cal Q}}(\varphi,\theta),\varphi,\theta,\ep)+J\tilde{{\cal S}}(\theta)F_2(\Psi+\tilde{{\cal Q}}(\varphi,\theta),\varphi,\theta,\ep)
\]
and
\[
{\mathscr N}_2(\Psi,\varphi,\theta,\ep)=F_2(\Psi+\tilde{{\cal Q}}(\varphi,\theta),\varphi,\theta,\ep).
\]
The functions ${\mathscr N}_{1,2}$ satisfy the symmetry properties (\ref{symproplat}) since the functions $F_{1,2}$ satisfy these properties, and
\[
\tilde{{\cal Q}}(\varphi+\pi/2,\theta)=R_{\varphi+\pi/2}\tilde{{\cal S}}(\theta)=R_{\pi/2}R_{\varphi}\tilde{{\cal S}}(\theta)=-J\tilde{{\cal Q}}(\varphi,\theta).
\]
Finally, setting $\hat{\varphi}=\varphi-\int_0^{\theta}\,h_2(s)ds+j\theta$ in the above system and dropping the hats gives us (\ref{preav_int}), where $V\equiv \hat{B}_j$ and
\[
\begin{array}{c}
H_1(\Psi,\varphi,\theta,\ep)=R_{\int_0^{\theta}\,h_2(s)ds}R_{-j\theta}{\mathscr N}_1(\Psi,\varphi+\int_0^{\theta}\,h_2(s)ds-j\theta,\theta,\ep)\\[0.2in]
H_2(\Psi,\varphi,\theta,\ep)={\mathscr N}_2(\Psi,\varphi+\int_0^{\theta}\,h_2(s)ds-j\theta,\theta,\ep),
\end{array}
\]
which satisfy the symmetry properties (\ref{symproplat}).  This ends the proof of Proposition \ref{propint}.
\hfill\qed

\end{document}